\input ./style/arxiv-general.cfg
\documentclass[aos,MSNbibl,nameyear,dvips]{arximspdf}
\makeatletter
   \@ifpackageloaded{graphicx}{}{\usepackage{graphicx}}
\makeatother
\usepackage{multirow,dcolumn}
\usepackage{url,breakurl}

\doi{10.1214/15-AOS1349}
\volume{44}
\issue{1}
\pubyear{2016}
\firstpage{58}
\lastpage{86}
\docsubty{FLA}

\makeatletter
\newcolumntype{d}[1]{D{.}{.}{#1}}
\def\mid{|}
\newcommand{\eqref}[1]{(\ref{#1})}
\def\overset{\stackrel}
\newcommand{\N}{\mathbb{N}}
\newcommand{\R}{\mathbb{R}}
\renewcommand{\P}{\mathbb{P}}
\newcommand{\E}{\mathbb{E}}
\newcommand{\Var}{\operatorname{\mathbb{V}{ar}}}
\newcommand{\Cov}{\operatorname{\mathbb{C}{ov}}}

\newtheorem{theorem}{Theorem}
\newtheorem{lemma}{Lemma}
\newtheorem{corollary}{Corollary}
\newproclaim{definition}{Definition}
\newproclaim{assumption}{Assumption}

\newproclaim{example}{Example}

\makeatother

\begin{document}
\begin{frontmatter}

\title{Near-optimal estimation of jump activity in semimartingales}
\runtitle{Near-optimal estimation of jump activity}

\begin{aug}
\author[A]{\fnms{Adam D.}~\snm{Bull}\corref{}\ead[label=e1]{a.bull@statslab.cam.ac.uk}\thanksref{T1}}
\runauthor{A.~D. Bull}
\thankstext{T1}{Supported by EPSRC Grant EP/K000993/1.}
\affiliation{University of Cambridge}
\address[A]{Statistical Laboratory\\
University of Cambridge\\
Wilberforce Road\\
Cambridge CB3 0WB\\
United Kingdom\\
\printead{e1}}
\end{aug}

%
\received{\smonth{9} \syear{2014}}
%
\revised{\smonth{3} \syear{2015}}

%
\begin{abstract}
In quantitative finance, we often model asset prices as
semimartingales, with drift, diffusion and jump components. The jump
activity index measures the strength of the jumps at high
frequencies, and is of interest both in model selection and fitting,
and in volatility estimation. In this paper, we give a novel
estimate of the jump activity, together with corresponding
confidence intervals. Our estimate improves upon previous work,
achieving near-optimal rates of convergence, and good finite-sample
performance in Monte-Carlo experiments.
\end{abstract}

%
\end{frontmatter}

\section{Introduction}
\label{sec:int}

In quantitative finance, we often wish to model asset prices, for
example, to price options or evaluate investment strategies. Typically,
we assume that asset log-prices are given by a semimartingale; in
other words, the sum of drift, diffusion and jump processes. In the
following, we will be interested in the jump activity index, a
parameter which determines the strength of the jump process at high
frequencies.

The jump activity is important for two reasons. First, any
semimartingale model will make claims about the jump activity;
typically, the activity is either assumed known and fixed, or is a
free parameter to be estimated. Knowledge of the jump activity thus
informs our choice of model, and may allow us to fit it more
accurately.

Second, the jump activity controls the difficulty of estimating
another parameter of interest, the volatility. This parameter
measures the strength of the diffusion component of price movements,
and is often a key target for financial modellers. It is known that
under high jump activity, the volatility becomes harder to estimate;
this problem can be avoided using specialised volatility estimates,
but at the cost of making stronger assumptions.

Knowledge of the jump activity is thus important both for the analysis
of individual price records, to inform the choice of volatility
estimate; and more generally in research, to guide the development of
future estimates. In the following, we will therefore investigate the
problem of accurately estimating the jump activity.

Previous attempts to estimate the jump activity of semimartingales
have either achieved poor rates of convergence, or worked only under
restrictive assumptions. In this paper, we will describe a new jump
activity estimate, which achieves near-optimal rates of convergence in
a general setting, along with improved finite-sample performance. We
will further provide limiting distributions for our estimate,
validated by Monte-Carlo experiments.

We begin by discussing in more detail the nature of the problem, and
relevant work in the literature. We will suppose we have a log-price
process given by a semimartingale $X_t$ on $[0, 1]$, and make
$n$ observations
\[
X_{j/n},\qquad  j= 0, \ldots, n-1.
\]
We then define the
jump activity index
\[
\beta= \inf \bigl\{r \in[0,2]: S(r) < \infty\ \mathrm{a.s.} \bigr\},\qquad S(r) = \sum
_{s \in[0,1]} \vert\Delta X_s
\vert^r,
\]
letting $\Delta X_s = X_s - X_{s^-}$ denote the jumps of $X_t$, and
using the
convention $0^0 = 0$.

As semimartingales have finite quadratic variation, we have $S(2) <
\infty$ almost surely, and so the jump activity $\beta\in[0,2]$.
When the sample path of $X_t$ has finitely many jumps, $\beta=
0$; when it may have infinitely many jumps, but the jumps are of
finite variation, $\beta\in[0,1]$; and when the jumps may be of
infinite variation, $\beta\in[1, 2]$. The more activity $X_t$
has in its small jumps, the larger we will have to choose $r$ to
make $S(r)$ finite, and the larger $\beta$ will be.

From Lemma~3.2.1 of \citet{jacoddiscretization2011}, we can
equivalently define
\[
\beta= \inf \bigl\{r \in[0,2]: I(r) < \infty\ \mathrm{a.s.} \bigr\},\qquad I(r) = \int
_0^1 \int_\R1 \wedge\vert
x\vert^r \nu(dx,ds),
\]
letting
$\nu(dx,ds)$ denote the compensator of the jump measure of $X_t$.
When $X_t$ is a L\'evy process, $\beta$ is thus the
Blumenthal--Getoor index [\citet{blumenthalsample1961}]; for example,
if $X_t$ is a stable process, then $\beta$ is its stability
parameter. More generally, $\beta$ gives an extension of the
Blumenthal--Getoor index to semimartingales.

The jump activity $\beta$ is thus a parameter of interest when
choosing models for the log-price process $X_t$. Many common models
assume either that no jumps are present, or that there are
finitely-many jumps almost-surely; in either case, we therefore assume
that $\beta= 0$. This includes all It\^ o process models, as well as
the Merton, Kou and Bates models, for example.

Some models allow positive values of $\beta$; for example, the
(time-changed) normal-inverse Gaussian, Meixner and generalised
hyperbolic models assume $\beta= 1$, while the (time-changed) CGMY
or tempered-stable model includes $\beta$ as a free parameter to be
estimated.
Knowledge of $\beta$ thus allows us to better decide between
competing models, and in the latter case also to fit these models to
price data. [For definitions of the models, see
\citet{contfinancial2004,papapantoleonintroduction2008}.]

Further interest in the jump activity arises from the problem of
volatility estimation. Let $X_t^c$ denote the continuous part of
$X_t$. Then the integrated volatility of $X_t$ over $[0,1]$,
given by the quadratic variation $[X^c]_1$, is a parameter of much
interest in options pricing or risk modelling, and its estimation has
been extensively studied.

When $X_t$ is continuous, the integrated volatility can be estimated
by the observed quadratic variation; however, price data is widely
accepted to contain jumps, which must be accounted for explicitly.
Methods for doing so include thresholding
[\citeauthor{mancinidisentangling2001}
(\citeyear{mancinidisentangling2001,mancininon2006})], bipower variation
[\citet{barndorff-nielsenpower2004,barndorff-nielsenlimit2006}], and
characteristic functions
[\citeauthor{todorovinverse2012}
(\citeyear{todorovinverse2012,todorovrealized2012})].

Unfortunately, the convergence rates of these methods suffer when the
jumps are of infinite variation. While this can be avoided if we
assume the jumps are driven by a stable-like process
[\citet{jacodefficient2014}], or that prices are given by a
time-changed process [\citet{bullestimating2013}], it is known that in
general, poor rates are unavoidable [\citet{jacodremark2012}].

When estimating volatility, we would therefore like to know whether
the jumps are of infinite variation, and if so, how active they are;
equivalently, we would like to know whether $\beta$ is greater than
1, and if so, by how much. This question is of interest both when
choosing an estimator to apply to particular price data, and also more
generally when planning research on volatility estimation.

Previous authors have attempted to recover $\beta$ in a variety of
settings, including when no diffusion component is present
[\citet{zhaononparametric2009,todorovactivity2010,woerneranalyzing2011}],
or when testing if $\beta$ is greater than zero
[\citet{leedetecting2010,ait-sahaliatesting2012}] or one
[\citet{contnonparametric2011}]. In the following, however, we will
concentrate on estimating $\beta$ in general, when a diffusion term
may also be present.

In this context, \citet{ait-sahaliaestimating2009} provide an
estimate of $\beta$ based upon jump counting. While
\citeauthor{ait-sahaliaestimating2009} cannot prove results for all
semimartingales, they do provide convergence rates under the
additional assumption that the jumps are dominated by a stochastic
integral of a stable-like process. Similar assumptions have also been
considered by \citet{jingestimating2011} and
\citet{jacodefficient2014}, for example, and are satisfied by many
common models of price data.

Under these conditions, \citet{ait-sahaliaestimating2009} show that
their estimate of the jump activity $\beta$ converges at a rate
$n^{-\beta/10}$. Related estimates have been considered also by
\citet{jingestimating2011}, \citet{ait-sahaliaidentifying2012} and
\citet{jingjump2012}; the best convergence is obtained by the
latter, who achieve the rate $n^{-\beta/8}$. However, this still
falls short of the corresponding lower bound of
$n^{-\beta/4}\log(n)^{-(1-\beta/4)}$, given by
\citet{ait-sahaliaidentifying2012}.

If we assume not only stable-like jumps, but also that log-prices are
given by a L\'evy process, Rei{\ss} (\citeyear{reistesting2013}) shows we can
estimate $\beta$ at the near-optimal rate
$n^{-\beta/4+\varepsilon}$, for any $\varepsilon> 0$. However,
the assumption of L\'evy behaviour is quite restrictive in a financial
context, and unfortunately the approach of
\citeauthor{reistesting2013} does not easily generalise to
semimartingales.

In the following, we will therefore describe a new estimate of the
jump activity $\beta$, using a multi-scale jump-counting
approach. We will show that by combining jump-counting estimates
across different time-scales, we will be able to cancel out the bias
in these estimates, obtaining improved accuracy.

In a similar setting to that of \citet{ait-sahaliaestimating2009},
with no assumption of L\'evy behaviour, our estimate will obtain the
near-optimal convergence rates $n^{-\beta/4+\varepsilon}$, as well
as improved finite-sample performance. We will also give limiting
distributions, validated by Monte-Carlo experiments.

In Section~\ref{sec:est}, we describe our estimates in full, and in
Section~\ref{sec:the}, discuss their theoretical properties. In
Section~\ref{sec:exp}, we then perform our Monte-Carlo experiments,
and in
Section~\ref{sec:pro}, give proofs.

\section{Jump activity estimates}
\label{sec:est}

We now describe our estimate of the jump activity $\beta$. In the
following, we will suppose that $\beta> 0$; we note the case
$\beta= 0$ can be tested for separately, for example, using the
methods of \citet{leedetecting2010} or
\citet{ait-sahaliatesting2012}.

Our approach builds upon the work of
\citet{ait-sahaliaestimating2009}, who estimate $\beta$ by
counting jumps in $X_t$. The authors define the jump counts
\[
\widetilde A_n(\tau) = \sum_{j=0}^{n-2}
1_{\tau
\vert X_{(j+1)/n}-X_{j/n}\vert \ge1},
\]
which for suitable $\tau> 0$,
approximate the number of jumps in $X_t$ of size at least
$\tau^{-1}$.

For $\rho> 1$, \citeauthor{ait-sahaliaestimating2009} then
estimate $\beta$ by
\[
\widehat\beta_n^{\mathrm{AJ}} = \log_{\rho} \biggl(
\frac{\widetilde A_n(\rho\tau_n)}{\widetilde A_n(\tau_n)} \biggr),
\]
using the convention $0/0=1$. If the jumps of $X_t$ are dominated
by a stochastic integral of a stable-like process, then as $n \to
\infty$, for suitable sequences $\tau_n$, we can expect
\[
\widetilde A_n(\tau_n) \approx C\tau_n^\beta,
\]
for some quantity $C > 0$. We would then have that
\[
\widehat\beta_n^{\mathrm{AJ}} \approx\log_\rho \biggl(
\frac{C(\rho\tau_n)^\beta}{C\tau_n^\beta} \biggr) = \beta.
\]

Unfortunately, \citeauthor{ait-sahaliaestimating2009} were not able
to provide good convergence rates for this method, as the estimates
$\widehat\beta^{\mathrm{AJ}}_n$ are too biased when $\tau_n$ is large. In
the following, we will therefore provide an improved version of this
method, which corrects for the bias in $\widehat\beta^{\mathrm{AJ}}_n$,
achieving near-optimal rates of convergence.

We will use three techniques to correct for this bias. First, will
we symmetrise the data, correcting for bias due to high-activity,
asymmetrically-distributed jumps. Second, we will smooth the jump
counts, correcting for bias due to the roughness of the indicator
function $1_{\vert x\vert \ge1}$. Finally, and most importantly, we will
eliminate the remaining bias by cancelling between estimates at
different time-scales.

We first describe a procedure to symmetrise the process $X_t$, as
given, for example, in \citet{jacodefficient2014}. For
$j= 0, \ldots, n-3$, we define random variables
\[
\Delta X_{j,n} = (X_{(j+2)/n} - X_{(j+1)/n}) -
(X_{(j+1)/n} - X_{j/n}).
\]
We note that when $X_t$ is a L\'evy process, the random variables
$\Delta X_{j,n}$ are symmetric, even if the increments of $X_t$
are not. More generally, we may think of the $\Delta X_{j,n}$ as
symmetrised increments of the process $X_t$, across time intervals
of length $2/n$.

In the following, we will wish to work with increments of $X_t$
across different time-scales simultaneously. For $k = 0,1,\ldots,$
$j= 0, \ldots, n-2k-1$, we therefore also define random variables
\[
\Delta X_{j,k,n} = \sum_{l=0}^{k-1}
\Delta X_{j+2l,n}.
\]
We can similarly consider the $\Delta X_{j,k,n}$ to be symmetrised
increments of $X_t$, now across time intervals of length $2k/n$.

Next, we will replace the indicator function $1_{\vert x\vert \ge1}$
with a smooth function, similarly to \citet{jingjump2012}. We will
use a smooth function $1 - K(x)$, where the kernel $K: \R\to
[0,1]$ is an even Schwartz function, equal to one in a neighbourhood
of the origin. For example, in our experimental results, we will
choose
\[
\label{eq:k} K(x) = %
\cases{1, &\quad $\vert x\vert \le1$,\vspace*{2pt}
\cr
\displaystyle\biggl(1+\exp \biggl(\frac{1}{2-\vert x\vert}-\frac{1}{\vert x\vert
-1} \biggr)
\biggr)^{-1},& \quad $1 \le\vert x\vert \le2$,\vspace*{2pt}
\cr
0, &\quad $\vert x
\vert \ge2$.} %
\]
We will also fix a constant $m \in\N$, giving the number of
time-scales to use for bias correction.

For $\tau> 0$, we then define the jump counts
\[
\widehat A_n(\tau) = 0 \vee\widehat A_n'(
\tau),\qquad \widehat A_n'(\tau) = \sum
_{j=0}^{n-2m-1} \widehat a_{j,n}( \tau),
\]
where for $j= 0, \ldots, n-2m-1$, we set
\[
\widehat a_{j,n}(\tau) = \sum_{k=1}^{m}
w_k \bigl(1-K(\tau\Delta X_{j,k,n}) \bigr),\qquad  w_k
= \frac{(-1)^{k+1}}{2k}\pmatrix{m \cr k}.
\]
For a constant $\rho> 1$, and
sequence $\tau_n > 0$, we finally estimate $\beta$ by
\[
\widehat\beta_n = 0 \vee\log_{\rho} \biggl(
\frac{\widehat
A_n(\rho\tau
_n)}{\widehat
A_n(\tau_n)} \biggr) \wedge2,
\]
using the
convention $0/0=1$.

When $m = 1$, this estimate is similar to the jump-counting estimate
of \citet{ait-sahaliaestimating2009}: we replace the increments
$X_{(j+1)/n}-X_{j/n}$ with symmetrised increments $\Delta
X_{j,n}$; replace the indicator function $1_{\vert x\vert \ge1}$ with a
smooth function $1-K(x)$; and clip the estimate $\widehat\beta_n$
to the interval $[0,2]$. When $m > 1$, we additionally replace
$\widehat A_n(\tau)$ with a linear combination of jump counts across
different time-scales, clipped to be nonnegative.

We note the clipping of $\widehat A_n(\tau)$ and $\widehat
\beta_n$ ensures that the estimate $\widehat\beta_n$ is always
reasonable, even when the jump counts $A_n'(\tau)$ may be
inaccurate. While this step makes no contribution to the asymptotic
behaviour of $\widehat\beta_n$, it does reduce its error in finite
time.

In the following sections, we will show that our changes reduce the
bias in the estimate $\widehat\beta_n$, providing both theoretical
and experimental improvements to accuracy. We will also use these
results to motivate the selection of parameters in our estimate: the
number of time-scales $m$, inverse jump threshold $\tau_n$, and
threshold ratio $\rho$.

We will further give limiting distributions for $\widehat\beta_n$,
allowing us to build confidence intervals for $\beta$. Define the
constants
\begin{eqnarray*}
C_{\beta,\rho} &=& K_{\beta,\rho}/\rho^\beta\log(
\rho)^2K_\beta ^2,
\\
K_{\beta} &=& \int_\R \bigl(1- K(x) \bigr)\vert x
\vert^{-(1+\beta)} \,dx,
\\
K_{\beta,\rho} &=&\int_\R \bigl(K(x)-K(\rho x)
\bigr)^2\vert x\vert^{-(1+\beta)} \,dx,
\end{eqnarray*}
and for $\widehat\beta_n \in(0, 2]$, the random variables
\[
\widehat U_n(\beta) = \frac{\tau_n^{\widehat\beta_n/2}(\widehat
\beta
_n -
\beta)}{\widehat\sigma_{\rho,n}},\qquad \widehat
\sigma_{\rho,n}^2 = \frac{C_{\widehat\beta_n,\rho}K_{\widehat\beta_n}\tau
_n^{\widehat
\beta_n}}{\widehat A_n(\tau_n)}.
\]
When $\widehat\beta_n = 0$, likewise define
\[
\widehat U_n(\beta) = -\infty.
\]
We note that the random variables $\widehat U_n(\beta)$ are always
well defined, as $\widehat\beta_n$ must lie within $[0,2]$.

We will be able to show that, under suitable conditions, the
standardised errors
\[
\widehat U_n(\beta) \overset{d}\to N(0, 1).
\]
We will therefore be able to define $\gamma$-level confidence
intervals for $\beta$,
\[
\widehat I_n(\gamma) = \bigl\{ \beta\in(0,2): \bigl\vert\widehat
U_n(\beta )\bigr\vert \le \Phi^{-1} \bigl(\tfrac{1}2(1+
\gamma) \bigr) \bigr\},
\]
where $\Phi$ denotes the standard Gaussian distribution function.

We note that the integrals $K_{\beta}$ and $K_{\beta,\rho}$ can
usually be computed numerically. In the case where $\widehat
\beta_n$ is very small but nonzero, the integration of $K_{\beta}$
may be slow to converge, and it may be preferable to instead take
$\widehat\beta_n = 0$. In our experimental tests, we did so for
$\widehat\beta_n < 10^{-3}$.

\section{Theoretical results}
\label{sec:the}

To describe our theoretical results, we must first state our
assumptions. The assumptions will be very similar to those made by
\citet{jacodefficient2014}, and essentially require that the jumps
of the log-price process $X_t$ are dominated by a stochastic
integral of a stable-like process. Similar assumptions have also been
made by \citet{ait-sahaliaestimating2009} and
\citet{jingestimating2011}, and are satisfied by many common models
of price data; we refer to \citet{jacoddiscretization2011} for
definitions and notation.

\begin{assumption}
\label{ass:model}
We first assume we have a probability space $(\Omega, \mathcal F,
\mathbb P)$, with filtration $\mathcal F_t$, and a jump activity
index $\beta\in(0, 2)$. We then assume the log-price process
\[
X_t = \int_0^t b_s \,ds
+ \int_0^t c_s \,dB_s +
\int_0^t \gamma_s^+
\,dL^+_s + \int_0^t
\gamma_s^- \,dL^-_s + \int_0^t
\int_\R\delta_s(x) \mu(dx, ds),
\]
where:
\begin{longlist}[(iii)]
\item[(i)]$B_t$ is an adapted Brownian motion;
\item[(ii)] the adapted Poisson random measure $\mu(dx,ds)$ has intensity
$dx \,ds$, and is independent of $B_t$;
\item[(iii)] the L\'evy processes
\[
L_t^\pm= \int_0^t \int
_\R\delta^\pm(x) \bigl(\mu(dx,ds) -
1_{\delta^\pm(x) < 1} \,dx \,ds \bigr),
\]
for disjointly-supported functions
$\delta^{\pm}(x) \ge0$, $\int_\R1 \wedge\delta^\pm(x)^2 \,dx <
\infty$;
\item[(iv)] the predictable processes $b_s$ and $\gamma_s^\pm$ are
locally bounded; and
\item[(v)] the predictable function $\delta_s(x)$ has $\int_\R
1 \wedge\vert\delta_s(x)\vert^{\upsilon_1} \,dx$ locally bounded,
for a parameter
$\upsilon_1 \in(0, \beta/2)$.
\end{longlist}

We additionally assume the volatility process
\begin{eqnarray*}
c_t &= &
c_0 +
\int_0^{t} b^c_s \,ds +
\int_0^{t} H_s \,dB_s +
\int_0^{t} H'_s
\,dB_s'
\\
&&{}+ \int_0^{t} \int_\R
\delta^c_s(x) \bigl(\mu(dx,ds) - 1_{\vert\delta^c_s(x)\vert < 1} \,dx \,ds
\bigr),
\end{eqnarray*}
where:
\begin{longlist}[(iii)]
\item[(i)] the adapted Brownian motion $B'_t$ is independent of $B_t$
and $\mu(dx,ds)$;
\item[(ii)] the predictable processes $b^c_s$, $H_s$ and $H'_s$ are
locally bounded; and
\item[(iii)] the predictable function $\delta^c_s(x)$ has $\int_\R1
\wedge\delta^c_s(x)^2 \,dx$ locally bounded.
\end{longlist}

We next assume the processes $L_t^\pm$ are close to
one-sided $\beta$-stable processes. Let
\[
F^\pm(U) = \int_{\delta^\pm(x) \in U} \,dx
\]
denote the L\'evy
measures of the processes $L_t^\pm$, and for $x > 0$, let
\[
\overline F^\pm(x) = F^\pm \bigl((x, \infty) \bigr)
\]
denote their upper L\'evy distribution functions. We then require that
for $x \in(0,1)$,
\[
\bigl\vert\overline F^\pm(x) - \beta^{-1}x^{-\beta}\bigr\vert =
O \bigl(x^{-\upsilon_2} \bigr),
\]
for a
parameter $\upsilon_2 < \beta- 1$.

Finally, we assume that the characteristics $b$, $H$ and
$\gamma^\pm$ are smooth in quadratic mean: we assume there are
stopping times $T_n \to\infty$, such that for $V = b$, $H$ or
$\gamma^\pm$, and any $0 \le t \le t+h \le1$,
\[
\E \bigl[(V_{(t+h) \wedge T_n} - V_{t \wedge
T_n})^2\mid\mathcal
F_{t} \bigr] = O(h),
\]
uniformly in $t$.
\end{assumption}

In other words, we assume that the log-price process $X_t$ and
volatility process $c_t$ are It\^o semimartingales; that the jumps
of $X_t$ are dominated by stochastic integrals against L\'evy
processes $L_t^\pm$, whose L\'evy distribution functions approach
those of a $\beta$-stable process; and that the drift process
$b_t$, leverage process $H_t$, and jump integrands
$\gamma^\pm_t$ exhibit smoothness behaviour typical of It\^o semimartingales.

We note that the jump processes in our assumptions are all described
using a Grigelionis representation, as integrals against a common
Poisson random measure~$\mu$; however, this condition is not
restrictive, as any collection of jump processes can be expressed in
this form [\citet{jacoddiscretization2011}, Theorem~2.1.2]. We
likewise note that while our assumptions choose a specific
normalisation for the jump processes $L^\pm_t$, this is not
restrictive, as the processes can always be rescaled by the terms
$\gamma^\pm_t$.

While the driving L\'evy processes $L_t^\pm$ must have stable-like
behaviour, our model allows for deviations from stability both in the
L\'evy distribution functions $\overline F^\pm$, which must be close
to stable only for small jumps; and in the idiosyncratic jumps
described by $\delta_t(x)$, which can account for any additional
jump activity. The presence of two separate one-sided L\'evy processes
$L_t^\pm$ also allows us to describe processes with asymmetric
jump activity.

We further allow the volatility $c_t$ to contain jumps and leverage,
and the other characteristic processes $b_t$, $H_t$ and
$\gamma_t^\pm$ to display a wide range of semimartingale
behaviour. Finally, we note that when the processes $\gamma^\pm_t$
are not both almost-surely zero, the parameter $\beta$ in our
assumptions agrees with the jump activity index as defined in
the \hyperref[sec:int]{Introduction}.

Under these assumptions, we will be able to provide limiting
distributions for the estimates $\widehat\beta_n$, and standardised
errors $\widehat U_n(\beta)$; we begin by defining the appropriate
notion of
convergence. Let $Z_n \in\R^d$ be random variables on a
probability space $(\Omega, \mathcal F, \mathbb P)$, and $Z \in
\R^d$ a random variable defined on a suitable extension $(\widetilde
\Omega, \widetilde{\mathcal F}, \widetilde{\mathbb P})$. We will say
the $Z_n$ converge stably in distribution to $Z$, $Z_n
\overset{\mathrm{sd}}\to Z$, if
\[
\E \bigl[Yf(Z_n) \bigr] \to\widetilde\E \bigl[Yf(Z) \bigr],
\]
for all random variables $Y \in\R$ on $\Omega$, and bounded
continuous functions $f:\R^d\to\R$
[\citet{jacoddiscretization2011}, Section~2.2.1].

We note that stable convergence in distribution is stronger than the
usual convergence in distribution, and requires convergence to hold
even after conditioning on events in $\mathcal F$. Under suitable
conditions, this stronger notion of convergence will allow us to show
not only that the estimates $\widehat\beta_n$ converge to unbiased
Gaussian mixtures, but also that the standardised errors $\widehat
U_n(\beta)$ converge to standard Gaussians.

To be precise, we first define the jump activity processes
\[
\overline\Gamma_t = \int_0^t
\overline\gamma_s \,ds,\qquad \overline \gamma_t =
\tfrac{1}2 \bigl(\bigl\vert\gamma^+_{t}\bigr\vert^\beta+ \bigl\vert
\gamma ^-_{t}\bigr\vert^\beta \bigr);
\]
we note that the process $\overline\gamma_t$ measures the
instantaneous stable-like jump activity at time $t$, while
$\overline\Gamma_1$ counts the total stable-like jump activity over
the interval $[0,1]$. We then have the following results.

\begin{theorem}
\label{thm:lim}
Under Assumption~\ref{ass:model}, let $m \in\N$, $\alpha= m/2(m+1)$,
$\tau_n = Cn^\alpha$ for some $C > 0$, and $\rho> 0$. Then on
the event $\overline\Gamma_1 > 0$,
\[
\tau_n^{\beta/2}(\widehat\beta_n - \beta)
\overset{\mathrm{sd}}\to\sigma_{\beta,\rho}Z, \qquad\widehat \sigma^2_{\rho, n}
\overset{p}\to\sigma_{\beta,\rho}^2, \qquad \widehat U_n(
\beta) \overset{\mathrm{sd}}\to Z,
\]
where the variance
\[
\sigma_{\beta,\rho}^2 = C_{\beta,\rho}/\overline
\Gamma_1,
\]
and $Z$ is a random variable
defined, on a suitable extension of the probability space $(\Omega,
\mathcal F, \mathbb P)$, to be standard Gaussian given $\mathcal
F$.
\end{theorem}

\begin{corollary}
\label{cor:lim}
In the setting of Theorem~\ref{thm:lim}, let $\gamma\in(0, 1)$. If
$\overline\Gamma_1$ is not almost surely zero, then
\[
\P \bigl(\beta\in\widehat I_n(\gamma) \mid\overline
\Gamma_1 > 0 \bigr) \to\gamma,
\]
and on the event $\overline\Gamma_1 > 0$, $\widehat
I_n(\gamma)$ has diameter $O_p(\tau_n^{-\beta/2})$.
\end{corollary}

We conclude that on the event that $X_t$ includes any stable-like
jump activity, the estimate $\widehat\beta_n$ converges at a
rate
\[
n^{-\beta m/4(m+1)};
\]
of course, we cannot expect convergence when no stable-like jumps are
present. Furthermore, on this event the $\widehat I_n(\gamma)$ are
indeed $\gamma$-level confidence intervals for $\beta$,
contracting at the correct rate.

A single-scale procedure, with $m = 1$, can thus converge at a rate
$n^{-\beta/8}$, recovering the results of
\citet{jingjump2012}. Moreover, by choosing $m$ large enough,
a~multi-scale procedure can achieve a rate $n^{-\beta/4+\varepsilon}$,
for any $\varepsilon> 0$. Indeed, this rate is near-optimal: a
corresponding lower bound rate of
$n^{-\beta/4}\log(n)^{-(1-\beta/4)}$ is given by
\citet{ait-sahaliaidentifying2012}.

Asymptotically, the rate of convergence will always be improved by
choosing $m$ larger, and we should therefore choose $m$ as large
as possible. In finite time, however, a larger choice of $m$ may
take longer to reach the asymptotic regime, and so be less accurate in
practice.

The optimal finite-time choice of $m$ may be difficult to compute,
and depends on a number of unknown quantities. However, in our
Monte-Carlo experiments, we found the choice $m = 3$ performed well,
and should already provide improvements over a single-scale
estimate. If more accuracy is desired, users may wish to perform a
simulation study to select $m$, or compare estimates for a number of
different choices of $m$.

Another practical consideration comes from microstructure noise. When
observing price data at high frequencies, it is now widely accepted
that observations of efficient prices are corrupted by
noise. Similarly to \citet{ait-sahaliaestimating2009}, for
reasonable choices of $\tau_n$ the estimates $\widehat\beta_n$
are sensitive only to large jumps in prices, and so will not be much
affected by noise. However, modifications to account for noise
asymptotically are also possible, as in \citet{jingestimating2011}
or \citet{bullestimating2013}, and may be left for future work.

\section{Monte-Carlo experiments}
\label{sec:exp}

We now perform Monte-Carlo tests of our multi-scale estimates
$\widehat\beta_n$, comparing them to the jump-counting estimates
$\widehat\beta_n^{\mathrm{AJ}}$ of \citet{ait-sahaliaestimating2009}. We
note that as $\widehat\beta_n^{\mathrm{AJ}}$ can sometimes be very large,
its RMSE can be distorted by the small chance of a large error. To
provide a fair comparison, we will therefore consider the clipped
estimates
\[
\widetilde\beta_n = 0 \vee\widehat\beta_n^{\mathrm{AJ}}
\wedge2,
\]
defined similarly to $\widehat\beta_n$; we note that this clipping
can only reduce the error in $\widehat\beta_n^{\mathrm{AJ}}$.

We will also compare our confidence intervals $\widehat I_n(\gamma)$
to similar ones defined in terms of $\widetilde\beta_n$. From
Theorem~3 of \citet{ait-sahaliaestimating2009}, and arguing as in
our Theorem~\ref{thm:lim}, we have that the
\[
\widetilde I_n(\gamma) = \bigl\{\beta\in(0,2): \bigl\vert\widetilde
U_n(\beta )\bigr\vert \le \Phi^{-1} \bigl(\tfrac{1}2(1+
\gamma) \bigr) \bigr\}
\]
are $\gamma$-level confidence intervals for $\beta$,
where
\[
\widetilde U_n(\beta) = \log(\rho) \bigl(\widetilde
A_n(\tau_n)^{-1} - \widetilde
A_n(\rho\tau_n)^{-1} \bigr)^{-1/2} (
\widetilde\beta_n - \beta).
\]
We may thus compare the $\widehat I_n(\gamma)$ to the $\widetilde
I_n(\gamma)$.


In each run of our simulation, we will generate $n = 23\mbox{,} 400$
observations, corresponding to observations taken every second of a
typical 6.5-hour trading day. Our observations will be drawn from a
log-price process
\[
X_t = B_t + \gamma_t R_t, \qquad t
\in[0,1],
\]
where $B_t$ is a standard Brownian motion; the deterministic scaling
process
\[
\gamma_t = (2t-1) \vee0;
\]
and the jump process
\[
R_t = \theta_1 S_t^\beta+
\theta_2 S_t^{\beta-0.2},
\]
for constants $\theta_1, \theta_2 > 0$, and independent
$\alpha$-stable processes $S_t^\alpha$.

The process $X_t$ thus models a price process with both diffusion
and jump components. Its jumps are driven by a $\beta$-stable
process, with time-varying intensity $\gamma_t$, but also contain a
nuisance component, with jump activity $\beta- 0.2$.

The relative strengths of these jumps are given by the constants
$\theta_1$ and $\theta_2$, which we will set in terms of a
parameter $p \in(0,1)$. To set $\theta_1$, we will require that
an increment $\theta_1(S_{j/n}^\beta- S_{(j-1)/n}^\beta)$ contains
a jump larger than 0.2 with probability $p$. To set $\theta_2$, we
will likewise require this condition holds for $\theta_2
S_t^{\beta-0.2}$, with probability $0.05p$.

To model the microstructure noise present at one-second time scales,
we will generate observations
\[
Z_j = X_{j/n} + 0.01\varepsilon_j,\qquad j = 0,
\ldots, n-1,
\]
where the independent noises $\varepsilon_j \sim N(0, 1)$. As noted
in \citet{ait-sahaliaestimating2009} and in Section~\ref{sec:the}, the
estimates $\widetilde\beta_n$ and $\widehat\beta_n$ can be
expected to be robust to the presence of such noise, and we will
compute them as though the observations $Z_j$ were noiseless.

The estimates $\widetilde\beta_n$ and $\widehat\beta_n$ then
depend on a number of parameters; we begin by considering the inverse
thresholds $\tau_n = Cn^\alpha$. As noted in
\citet{ait-sahaliaestimating2009}, $\tau_n$ should be chosen to
ensure our jump counts $\widetilde A_n$ or $\widehat A_n$ will be
zero when no jump is present; the constant $C$ should thus be chosen
relative to the size of the diffusion component of $X_t$, as
measured for example by its integrated volatility.

In our simulations, we know that the integrated volatility of $X_t$
is equal to one, and so we may choose our parameters accordingly. In
general, the volatility will not be equal to one; however, we can
achieve a similar effect by first renormalising the observations to
have estimated integrated volatility equal to one. Such an estimate
could be provided by the method of \citet{podolskijestimation2009},
for example, although we will not pursue this further here.

In any case, we may now choose our parameters without worrying about
issues of scale. With $\widetilde\beta$,
\citet{ait-sahaliaestimating2009} recommend a threshold rate
$\alpha= 1/5$; with $\widehat\beta$, we will instead use the
rate $\alpha$ given by Theorem~\ref{thm:lim}. It remains to choose the
constants $C$, $\rho$, and for $\widehat\beta$ also $m$; in
our tests, we found the values $C = 0.05$, $\rho= 2$, and $m =
3$ worked well.

\begin{table}
\caption{Simulated means and standard deviations
of the estimates $\widetilde\beta_n$ and
$\widehat\beta_n$, as well as coverages of the
95\% confidence intervals $\widetilde I_n(0.95)$
and $\widehat I_n(0.95)$}
\label{tab:mv}
\begin{tabular*}{\textwidth}{@{\extracolsep{\fill}}lccccccc@{}}
\hline
&& \multicolumn{3}{c}{$\bolds{\widetilde\beta}$} &
\multicolumn{3}{c@{}}{$\bolds{\widehat\beta}$}\\[-6pt]
&& \multicolumn{3}{c}{\hrulefill} &
\multicolumn{3}{c@{}}{\hrulefill}\\
$\bolds{\beta}$ &$\bolds{p}$ & \textbf{0.5\%} & \textbf{1\%} &
\textbf{2\%} & \textbf{0.5\%} & \textbf{1\%} & \textbf{2\%}
\\
\hline
0.4 & mean& 0.39 & 0.39 & 0.39& 0.39 & 0.39 & 0.39\\
& std. dev.& 0.18 & 0.13 & 0.09& 0.13 & 0.09 & 0.07\\
& 95\% cov.& 0.90 & 0.91 & 0.93& 0.91 & 0.92 & 0.91\\[3pt]
0.8 & mean& 0.85 & 0.82 & 0.80& 0.81
& 0.80 & 0.79\\
& std. dev.& 0.49 & 0.34 & 0.23& 0.26 & 0.18 & 0.13\\
& 95\% cov.& 0.88 & 0.92 & 0.93& 0.92 & 0.93 & 0.93\\[3pt]
1.2 & mean& 1.13 & 1.21 & 1.23& 1.22
& 1.22 & 1.20\\
& std. dev.& 0.75 & 0.59 & 0.46& 0.40 & 0.29 & 0.21\\
& 95\% cov.& 0.80 & 0.90 & 0.93& 0.93 & 0.94 & 0.92\\[3pt]
1.6 & mean& 0.91 & 1.26 & 1.44& 1.54
& 1.58 & 1.57\\
& std. dev.& 0.92 & 0.80 & 0.63& 0.43 & 0.36 & 0.31\\
& 95\% cov.& 0.53 & 0.77 & 0.89& 0.93 & 0.93 & 0.92\\
\hline
\end{tabular*}
\end{table}

Table~\ref{tab:mv} then gives the mean and standard deviation of 10,000
simulated estimates $\widetilde\beta_n$ or $\widehat\beta_n$,
for a number of choices of $\beta$, $p$ and $m$. The table also
gives the simulated coverage of the 95\% confidence intervals
$\widetilde I_n(0.95)$ or $\widehat I_n(0.95)$. We see that the
multi-scale estimate $\widehat\beta_n$ has reduced bias and
variance compared with the single-scale estimate $\widetilde
\beta_n$, while the confidence intervals $\widehat I_n(0.95)$ retain
good coverage, improving upon $\widetilde I_n(0.95)$ when $\beta$
is large.

\begin{figure}[b]

\includegraphics{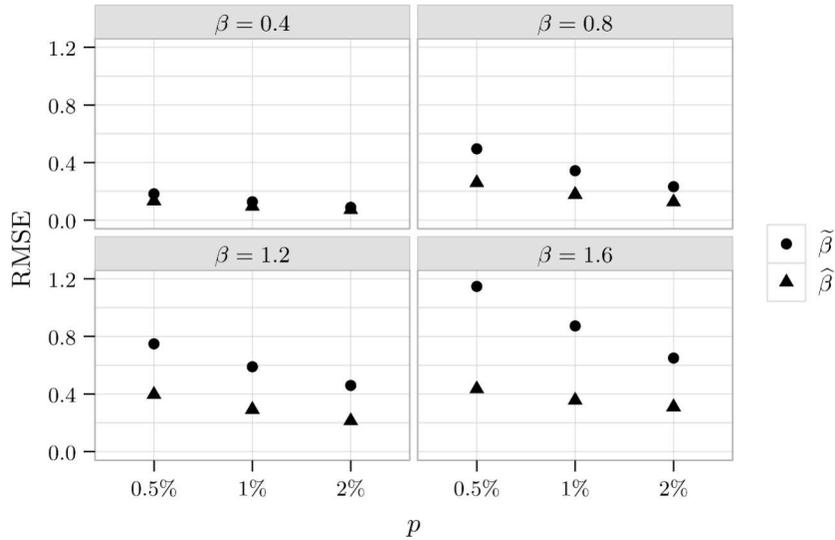}

\caption{Simulated RMSEs of the estimates $\widetilde\beta_n$ and
$\widehat\beta_n$.}
\label{fig:rmse}
\end{figure}

\begin{figure}

\includegraphics{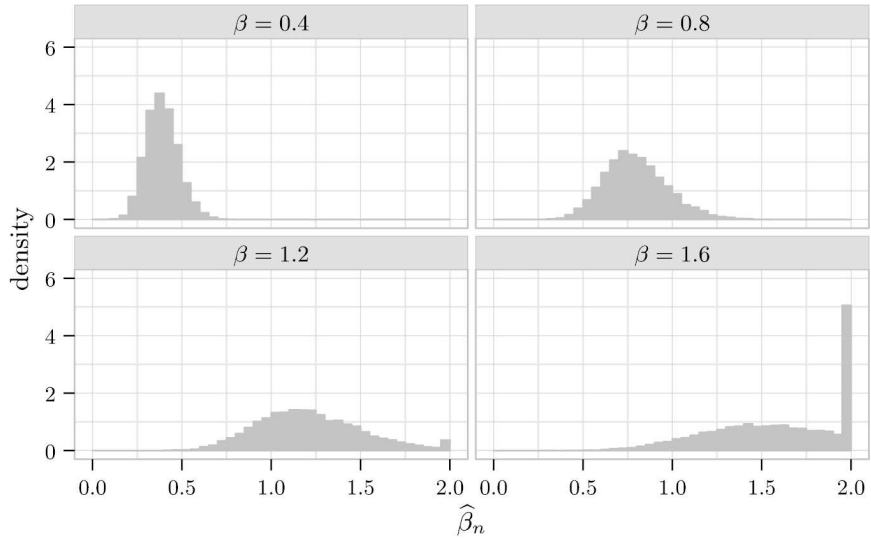}

\caption{Simulated distributions of the estimates $\widehat
\beta_n$, $p = 1\%$.}
\label{fig:beta}
\end{figure}

Figure~\ref{fig:rmse} plots the RMSE of the estimates $\widetilde
\beta_n$ and $\widehat\beta_n$; in the case $p = 1\%$,
Figure~\ref{fig:beta} further gives the full simulated distribution of
$\widehat\beta_n$. Again, we can see the multi-scale estimate
$\widehat\beta_n$ is more accurate than the single-scale estimate
$\widetilde\beta_n$. While the accuracy of $\widehat\beta_n$
suffers when $\beta$ is large, it remains good enough to distinguish
between different values of $\beta$.

\begin{figure}

\includegraphics{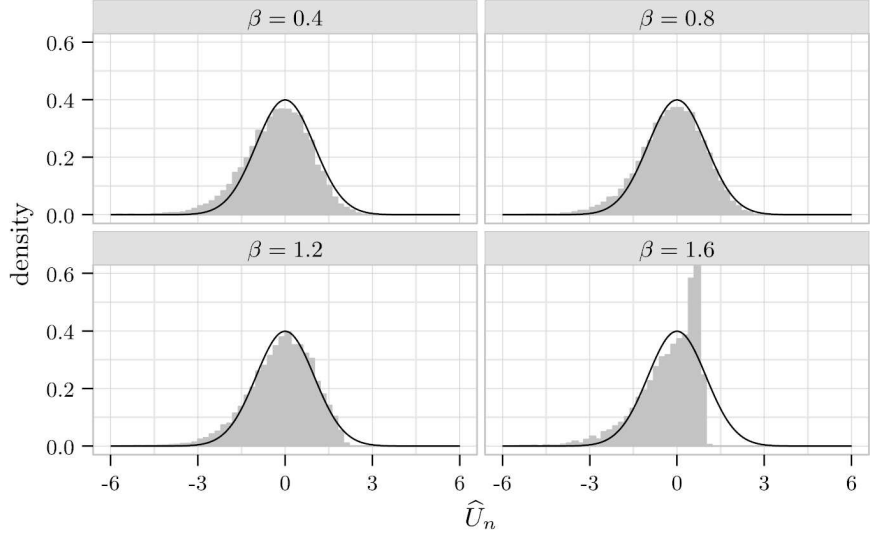}

\caption{Simulated distributions of the standardised errors $\widehat
U_n$, $p = 1\%$.}
\label{fig:z}
\end{figure}

Finally, Figure~\ref{fig:z} plots the simulated distribution of the
standardised errors $\widehat U_n$, together with the density of a
standard Gaussian distribution, shown as a solid line. We can see
that even in the finite-sample case, for $\beta= 0.4, 0.8, 1.2$,
the errors $\widehat U_n$ show good agreement with their asymptotic
distributions.

In the case $\beta= 1.6$, we see a strong deviation from Gaussian
on the right tail of $\widehat U_n$, due to the clipping of
$\widehat\beta_n$ at 2. This clipping, however, serves only to
reduce the error in the estimate $\widehat\beta_n$, and so does not
harm the coverage of the confidence intervals $\widehat
I_n(\gamma)$. Furthermore, the effect can be expected to disappear as
$n$ tends to infinity.

\section{Proofs}
\label{sec:pro}

We now give a proof of Theorem~\ref{thm:lim}. In Section~\ref
{sec:pro-tr}, we
will state the technical results we require; in Section~\ref{sec:pro-mp},
prove our main results; and in the supplementary material
[\citet{bullsupplement2015}], give the remaining technical
proofs.

\subsection{Technical results}
\label{sec:pro-tr}

We begin with a technical lemma bounding various stochastic integrals,
similarly to \citet{jacoddiscretization2011}.

\begin{lemma}
\label{lem:bounds}
Let $B_s$ be a Brownian motion, $\mu(dx,ds)$ a
Poisson jump measure with intensity $dx \,ds$, $a_s$ a
predictable process, $f_s(x)$ a predictable function, $t \in
[0,1]$, and $\kappa_p > 0$ denote constants depending only on $p
\ge1$.
\begin{longlist}[(iii)]
\item[(i)]
If $\int_0^t\vert a_s\vert \,ds < \infty$,
\[
\biggl\vert\int_0^t a_s \,ds
\biggr\vert^p \le t^{p-1} \int_0^t
\vert a_s\vert^p \,ds.
\]
\item[(ii)]
If $a_s$ is locally bounded,
\[
\E \biggl[\biggl\vert\int_0^t a_s
\,dB_s\biggr\vert^p \biggr] \le\kappa_p \E \biggl[
\biggl(\int_0^t a_s^2
\,ds \biggr)^{p/2} \biggr].
\]
\item[(iii)]
If $\int_\R f_s(x)^2 \,dx$ is locally
bounded, and $p \in[1,2]$, then
\[
\E \biggl[ \biggl\vert\int_0^t \int
_\R f_s(x) \bigl(\mu(dx,ds) - \,dx \,ds \bigr)
\biggr\vert^p \biggr] \le\kappa_p \E \biggl[\int
_0^t \int_\R \bigl\vert
f_s(x)\bigr\vert^p \,dx \,ds \biggr].
\]
\item[(iv)]
If $\int_\R1\wedge\vert f_s(x)\vert \,dx$
is locally bounded, then
\[
\E \biggl[1 \wedge\biggl\vert\int_0^t \int
_\R f_s(x) \mu(dx,ds)\biggr\vert^p
\biggr]\le\kappa_p\E \biggl[\int_0^t
\int_\R1 \wedge \bigl\vert f_s(x)\bigr\vert \,dx \,ds
\biggr].
\]
\item[(v)]
If $a_s$ and $\int_\R1\wedge
f_s(x)^2 \,dx$ are locally bounded, $p \in[1,2]$,
and $\alpha\ge0$, then
\begin{eqnarray*}
&&\E \biggl[1 \wedge\biggl\vert t^{-\alpha}\int_0^t
a_s\int_\R f_s(x) \bigl(\mu
(dx,ds) - 1_{\vert f_s(x)\vert < 1} \,dx \,ds \bigr)\biggr\vert^p \biggr]
\\
&&\qquad\le\kappa_p\E \biggl[ %
\int_0^t \bigl\vert t^{-\alpha}a_s
\bigr\vert ^p \int_{\vert f_s(x)\vert < t^\alpha} \bigl\vert f_s(x)
\bigr\vert^p \,dx \,ds +\int_0^t \int
_{\vert f_s(x)\vert \ge t^\alpha} \,dx \,ds
\\
&&\qquad\quad{}+ t^{p-1}\int_0^t\biggl \vert
t^{-\alpha}a_s \int_{t^{\alpha} \le\vert f_s(x)\vert < 1}
f_s(x) \,dx\biggr\vert^p \,ds \biggr].
\end{eqnarray*}
\end{longlist}
\end{lemma}

\begin{pf}
Parts~(i) and~(ii) are immediate from the
H\"older and Burkholder--Davis--Gundy inequalities,
respectively. Part~(iii) follows from Lemma~2.1.5 of
\citet{jacoddiscretization2011}, and part~(iv) likewise
follows from their Lem\-ma~2.1.8, noting that the left-hand side is
decreasing in $p$.

Finally, let $W$ denote the left-hand side of
part~(v). We make the decomposition
\begin{eqnarray*}
W &\le&\kappa_p\E \biggl[ %
\biggl\vert\int_0^t \int
_\R g_{1,s}(x) \bigl(\mu(dx,ds) - dx \,ds \bigr)
\biggr\vert^p
\\
&&{}+ 1 \wedge\biggl\vert\int_0^t \int
_\R g_{2,s}(x) \mu(dx,ds)\biggr\vert^p +
\biggl\vert\int_0^t \int_\R
g_{3,s}(x) \,dx \,ds\biggr\vert^p \biggr],
\end{eqnarray*}
where the terms
\[
g_{i,s}(x) = t^{-\alpha} a_s f_s(x)1_{J_i}
\bigl(\bigl\vert f_s(x)\bigr\vert \bigr),
\]
for intervals
\[
J_1 = \bigl[0, t^\alpha\bigr), \qquad J_2 = \bigl[t^\alpha,
\infty\bigr),\qquad  J_3 = \bigl[t^\alpha, 1\bigr).
\]
We deduce that
\begin{eqnarray*}
W &\le&\kappa_p \E \biggl[ %
\int_0^t \int_\R
\bigl\vert g_{1,s}(x)\bigr\vert^p \,dx \,ds +\int_0^t
\int_\R1 \wedge\bigl\vert g_{2,s}(x)\bigr\vert \,dx \,ds
\\
&&{}+ t^{p-1} \int_0^t\biggl \vert\int
_\R g_{3,s}(x) \,dx\biggr\vert^p \,ds \biggr],
\end{eqnarray*}
using parts~(i),~(iii) and~(iv). The
desired result follows.
\end{pf}

Next, we give a technical result on the characteristic exponents of
one-sided stable processes.

\begin{lemma}
\label{lem:ft}
Let $\mathcal F[f](u) = \int_{\R^d} \exp(i\langle
u,x\rangle)f(x) \,dx$ denote the Fourier transform, $\Gamma(x)$ the
gamma function, $\beta\in(0,2)$, $u \in\R$, and
\[
C_{\beta} = %
\cases{-2\Gamma(-\beta)\cos(\beta\pi/2), &\quad $
\beta\ne 1$,\vspace*{2pt}
\cr
\pi, &\quad $\beta= 1$.} %
\]
We then have:
\begin{longlist}[(ii)]
\item[(i)]
$\int_0^\infty(1 - \cos(ux))x^{-(1+\beta)} \,dx
= \tfrac{1}2 C_\beta\vert u\vert^\beta$; and
\item[(ii)]
$C_\beta\int_\R
\mathcal F[K](u) \vert u\vert^\beta \,du = 2\pi K_{\beta}.$
\end{longlist}
\end{lemma}

\begin{pf} We show each result in turn.
\begin{longlist}[(i)]
\item[(i)] This is a well-known result on stable processes; see, for
example,
Lem\-ma~14.11 of \citet{satolevy1999}.
%
\item[(ii)] For $\beta\ne1$, using generalised functions, we have
\begin{eqnarray*}
C_\beta\int_\R\mathcal F [K](u) \vert u
\vert^\beta \,du &=& C_\beta\int_\R K(x)
\mathcal F \bigl[\vert u\vert^\beta \bigr](x) \,dx
\\
&= &2\pi\int_\R \bigl(1-K(x) \bigr) \vert x
\vert^{-(1+\beta)} \,dx,
\end{eqnarray*}
since $K$ is symmetric, and $K(0) = 1$. For $\beta= 1$, the
same holds by analytic continuation.\quad\qed
\end{longlist}
\noqed\end{pf}

Using these lemmas, we will be able to prove several L\'evy
approximations to the behaviour of random variables
$\int_t^{t+h}a_s \,dX_s$. These approximations will hold under a
localisation assumption; by standard techniques,
we will be able to assume the following.

\begin{assumption}
\label{ass:bd} Assumption~\ref{ass:model} holds, the processes $b_t$,
$b^c_t$, $c_t$, $H_t$, $H'_t$, $\gamma^\pm_t$, $\int_\R
1 \wedge\vert\delta_t(x)\vert^{\upsilon_1} \,dx$ and $\int_\R1
\wedge\delta^c_t(x)^2 \,dx$ are uniformly bounded, and the stopping
time $T_1 = \infty$.
\end{assumption}

We now state our L\'evy approximation results; proofs of these results
will be given in the supplementary material
[\citet{bullsupplement2015}]. Our first result bounds the error in
approximating variables $\int_{t}^{t+h}a_s \,dX_s$ by L\'evy
integrals.

\begin{lemma}
\label{lem:lvy}
Under Assumption~\ref{ass:bd}, let $0 \le t \le t+h \le1$, set
\[
\xi_{t+h} = \int_{t}^{t+h}
a_s \,dX_s,
\]
for a deterministic real-valued process $a_s$ satisfying
$\vert a_s\vert \le1$,
and define the L\'evy approximation
\[
\overline\xi_{t+h} = \int_{t}^{t+h}
a_s \bigl(b_{t} \,ds + c_{t} \,dB_s +
\gamma_{t}^+ \,dL^+_s + \gamma_{t}^-
\,dL^-_s \bigr).
\]
Then the
approximation error
\[
\xi_{t+h} - \overline\xi_{t+h} = Y_{1} +
Y_{2},
\]
where the random variable
\[
Y_{1} = \int_t^{t+h} a_s
\biggl( %
H_{t}
B_s + \int_{t}^sH_{r}'
\,dB_r'
+ \int_t^s \int_{\vert\delta^c_r(x)\vert < 1}
\delta^c_r(x) \bigl(\mu(dx,dr) - dx \,dr \bigr) \biggr)
\,dB_s,
\]
and for $\alpha
\in(0, \tfrac{1}2)$, $u = O(h^{-\alpha})$, and some $\varepsilon
> 0$, we have
\[
\E \bigl[\vert uY_{1}\vert^2 \mid \mathcal
F_{t} \bigr] = O \bigl(h^{1+\varepsilon} \bigr),\qquad \E\bigl[1 \wedge\vert
uY_{2}\vert \mid\mathcal F_{t}\bigr] = O \bigl(h^{1+\varepsilon-\alpha\beta/2}
\bigr),
\]
uniformly over $a_s$ and
$t$.
\end{lemma}

Next, we state a result on the characteristic functions of random
variables $\int_{t}^{t+h}a_s \,dX_s$. Our argument will follow
Lemmas 11 and 12 of \citet{jacodefficient2014}, although we give a
tighter bound than in those results.

\begin{lemma}
\label{lem:cfs}
In the setting of Lemma~\ref{lem:lvy}, suppose also that $\vert
a_s\vert =
1$, and $\int_{t}^{t+h} a_s \,ds = 0$. Then for some $\varepsilon
> 0$, we have
\[
\E \bigl[\cos(u\xi_{t+h})\mid\mathcal F_{t} \bigr]= \exp
\biggl(-\int_{t}^{t+h}\theta_{t}(a_su)
\,ds \biggr) + O \bigl(h^{1+\varepsilon-\alpha\beta/2} \bigr),
\]
uniformly over
$a_s$ and $t$, where
\[
\theta_t(u) = \tfrac{1}2(c_t u)^2 +
C_\beta\overline\gamma_t \vert u\vert^\beta.
\]
\end{lemma}

Our final technical result gives a large-jump approximation to
functions of integrals $\int_{t}^{t+h} a_s \,dX_s$.

\begin{lemma}
\label{lem:jmp}
In the setting of Lemma~\ref{lem:lvy}, suppose $\vert a_s\vert =
1$, let
$t' \in[t, t+h]$, and set $h' = t+h - t'$. Also let $f$ be a
bounded even function, constant in a neighbourhood of the origin,
whose derivative $f'$ is a Schwartz function. Then
\[
\E \bigl[f(u \xi_{t+h})\mid\mathcal F_{t'} \bigr] = f(u
\xi_{t'}) + h' \vert u\vert^\beta\overline
\gamma_{t'} \int_\R \bigl(f(x)-f(0) \bigr)\vert
x \vert^{-(1+\beta)} \,dx + Y,
\]
for a term $Y$ satisfying $\E[\vert Y\vert \mid\mathcal F_{t}] =
o(h^{1-\alpha\beta})$,
uniformly in $a_s, t$ and $t'$.
\end{lemma}

\subsection{Main proofs}
\label{sec:pro-mp}

We now prove our main results. In the following, we will use the
shorthand
\[
t_j = j/n, \qquad t_{j,k} = (j+2k)/n.
\]
Our next lemma then bounds the means of our jump counts $\widehat
a_{j,n}(\tau)$.

\begin{lemma}
\label{lem:em}
Under Assumption~\ref{ass:bd}, for $m$ and $\tau_n$ as in the
statement of Theorem~\ref{thm:lim}, we have
\[
\E \bigl[\widehat a_{j,n}(\tau_n) \mid\mathcal
F_{t_j} \bigr] = \tau_n^\beta K_{\beta}
n^{-1}\overline\gamma_{t_j} + o \bigl(n^{-(1-\alpha\beta/2)}
\bigr),
\]
uniformly in $j= 0, \ldots, n-2m-1$.
\end{lemma}

\begin{pf}
We can equivalently define the constants $w_k$ by
%
\begin{equation}
\label{eq:wk} w_k = \sum_{l=k\vee1}^m(-1)^{k+1}(2l)^{-1}\pmatrix{l
\cr k},
\end{equation}
letting the above also define a new constant $w_0$. We then
have
\[
2\pi\widehat a_{j,n}(\tau_n) = 2\pi\sum
_{k=0}^{m} w_k \bigl(1 - K(
\tau_n\Delta X_{j,k,n}) \bigr),
\]
since the summand vanishes for $k=0$,
\[
= -2\pi\sum
_{k=0}^{m} w_k K(\tau_n
\Delta X_{j,k,n}),
\]
since $\sum_{k=0}^m
w_k = -\sum_{l=1}^m
(2l)^{-1} (1-1)^l = 0$,
\[
= -\int_\R
\mathcal F[K](u)\sum_{k=0}^m
w_k \cos(u\tau_n\Delta X_{j,k,n}) \,du,
\]
by Fourier inversion,
\[
= - \int_{\vert u\vert \le n^{\varepsilon}} \mathcal
F[K](u)\sum_{k=0}^m w_k
\cos(u\tau_n\Delta X_{j,k,n}) \,du+ O \bigl(n^{-1}
\bigr),
\]
for any $\varepsilon> 0$, since $K$ is Schwartz.

For small enough $\varepsilon$, setting $\theta_{j,n}(u) =
n^{-1}\theta_{t_j}(\tau_n u)$, we deduce
\begin{eqnarray*}
&&2\pi\E \bigl[\widehat a_{j,n}(\tau_n) \mid\mathcal
F_{t_j} \bigr]
\\
&&\qquad= -\int_{\vert u\vert \le n^{\varepsilon}} \mathcal F[K](u) \sum
_{k=0}^{m} w_k\exp \bigl(-2k
\theta_{j,n}(u) \bigr) \,du+ o \bigl(n^{-(1-\alpha\beta/2)} \bigr),
\end{eqnarray*}
using
Lemma~\ref{lem:cfs},
\begin{eqnarray*}
&=& \int
_{\vert u\vert \le n^\varepsilon
} \mathcal F[K](u) \sum_{l=1}^{m}(2l)^{-1}
\bigl(1-\exp \bigl(-2\theta_{j,n}(u) \bigr) \bigr)^l \,du
\\
&&{}+ o \bigl(n^{-(1-\alpha\beta/2)} \bigr),
\end{eqnarray*}
from \eqref{eq:wk},
\[
= \int_{\vert u\vert
\le n^\varepsilon} \mathcal F[K](u)
\bigl(\theta_{j,n}(u) + O \bigl(\theta_{j,n}(u)^{m+1}
\bigr) \bigr) \,du + o \bigl(n^{-(1-\alpha\beta/2)} \bigr),
\]
considering the Taylor series of $\log(1-x)$,
\[
= \int
_{\R} \mathcal F[K](u) \theta_{j,n}(u) \,du + o
\bigl(n^{-(1-\alpha\beta/2)} \bigr),
\]
since $K$ is Schwartz, and for $\vert u\vert
\le n^{\varepsilon}$, $\theta_{j,n}(u) = O (n^{-(1-2(\alpha+\varepsilon))}
)$,
\[
= n^{-1}\tau_n^\beta\overline
\gamma_{t_j}C_\beta\int_\R\mathcal
F[K](u)\vert u\vert^\beta \,du + o \bigl(n^{-(1-\alpha\beta/2)} \bigr),
\]
since $K$ is constant in a region of the origin, and so $\mathcal F[K]$ is
orthogonal to polynomials vanishing at the origin,
\[
= 2\pi\tau_n^\beta
K_{\beta} n^{-1}\overline\gamma _{t_j} + o
\bigl(n^{-(1-\alpha\beta/2)} \bigr),
\]
using Lemma~\ref{lem:ft}(ii).
\end{pf}

We next prove a lemma giving the variance of terms like
$K(\tau_n\Delta X_{j,k,n})$. To begin, for $\beta\in(0,2)$,
$\rho> 0$, we define the constants
\[
\overline K_{\beta,\rho} = \rho^{-\beta/2}\int_\R
\bigl(1- K(x) \bigr) \bigl(1-K(\rho x) \bigr)\vert x\vert^{-(1+\beta)} \,dx.
\]
We then have the following result.

\begin{lemma}
\label{lem:ev}
Under Assumption~\ref{ass:bd}, for $m$ and $\tau_n$ as in the
statement of Theorem~\ref{thm:lim}, let $j, j' = 0, \ldots, n-2m-1$, and
$k, k'= 1, \ldots, m$. Also let $a_s, a_s'$ be deterministic
processes with $\vert a_s\vert = \vert a_s'\vert = 1$, let $\max(t_j,
t_{j'}) \le t \le t+h \le\min(t_{j,k}, t_{j',k'})$, and set
\begin{eqnarray*}
\xi_{t'} &= &\int_{t_j}^{t'}
a_s \,dX_s,\qquad \xi_{t'}' = \int
_{t_{j'}}^{t'} a_s'
\,dX_s,
\\
V &= &\E \bigl[K(\tau_n \xi_{t_{j,k}})\mid\mathcal
F_{t+h} \bigr],\qquad V' = \E \bigl[K \bigl(\rho
\tau_n \xi_{t_{j',k'}}' \bigr)\mid\mathcal
F_{t+h} \bigr].
\end{eqnarray*}
Then
\[
\Cov \bigl[V,V'\mid\mathcal F_{t} \bigr] = h
\tau_n^\beta\rho^{\beta/2}\overline K_{\beta,\rho}
\overline \gamma _{t}+ Y,
\]
for a term $Y$ satisfying $\E[\vert Y\vert \mid\mathcal
F_{\min(t_j, t_{j'})}] = o(n^{-(1-\alpha\beta)})$, uniformly.
\end{lemma}

\begin{pf}
In the following, let $Y$ denote any term satisfying
\[
\E\bigl[\vert Y\vert \mid\mathcal F_{\min(t_j, t_{j'})}\bigr] = o \bigl(n^{-(1-\alpha\beta)}
\bigr).
\]
Repeatedly applying Lemma~\ref{lem:jmp},
we have
\begin{eqnarray*}
&&\E \bigl[VV' \mid\mathcal F_{t} \bigr]
\\
&&\qquad=\E \bigl[ %
 \bigl(K(
\tau_n \xi_{t+h}) - (t_{j,k}-t-h)
\tau_n^\beta\overline\gamma_{t+h}K_{\beta}
\bigr)
\\
&&\qquad\quad{}\times \bigl(K \bigl(\rho\tau_n \xi_{t+h}'
\bigr) - (t_{j',k'}-t-h) (\rho\tau _n)^\beta
\overline\gamma_{t+h}K_{\beta} \bigr) \mid\mathcal
F_{t} \bigr] + Y
\\
&&\qquad= %
\E \bigl[K(\tau_n
\xi_{t+h})K \bigl(\rho\tau _n \xi_{t+h}'
\bigr) \mid\mathcal F_{t} \bigr]
\\
&&\qquad\quad{}- \bigl((t_{j,k}-t-h) + \rho^\beta(t_{j',k'}-t-h)
\bigr)\tau_n^\beta \overline\gamma_{t}K_{\beta}
+ Y
\\
&&\qquad= %
K(\tau_n
\xi_{t})K \bigl(\rho\tau_n\xi _{t}'
\bigr)
\\
&&\qquad\quad{}- h\tau_n^\beta\overline\gamma_{t} \int
_\R \bigl(1-K(x)K(\rho x) \bigr)\vert x\vert^{-(1+\beta)}
\,dx
\\
&&\qquad\quad{}- \bigl((t_{j,k}-t-h) + \rho^\beta(t_{j',k'}-t-h)
\bigr)\tau_n^\beta\overline\gamma_{t}K_{\beta
}
+ Y.
\end{eqnarray*}
Again applying Lemma~\ref{lem:jmp}, we deduce that
\begin{eqnarray*}
&&\Cov \bigl[V, V'\mid\mathcal F_{t} \bigr]
\\
&&\qquad= \E \bigl[VV' \mid\mathcal F_{t} \bigr] -\E[V\mid
F_{t}]\E \bigl[V'\mid\mathcal F_{t} \bigr]
\\
&&\qquad= h\tau_n^\beta\overline\gamma_{t} \biggl(
\bigl(1+\rho^\beta \bigr)K_\beta- \int_\R
\bigl(1-K(x)K(\rho x) \bigr)\vert x\vert^{-(1+\beta)} \,dx \biggr) +Y
\\
&&\qquad= h\tau_n^\beta\rho^{\beta/2}\overline
\gamma_{t}\overline K_{\beta
,\rho
} + Y.
\end{eqnarray*}
\upqed\end{pf}

Next, we prove a lemma bounding the covariation of terms
$K(\tau_n\Delta X_{j,k,n})$ with other martingales.

\begin{lemma}
\label{lem:cor}
Under Assumption~\ref{ass:bd}, for $m$ and $\tau_n$ as in the statement
of Theorem~\ref{thm:lim}, let $t \in[0,1]$, and $k = 1, \ldots, m$. Then
\[
\sum_{j=0}^{\lfloor nt
\rfloor-2m}\E \bigl[ \bigl(1-K(
\tau_n\Delta X_{j,k,n}) \bigr) (M_{t_{j,k}}-M_{t_j})
\mid\mathcal F_{t_j} \bigr]= o_p \bigl(n^{\alpha\beta/2}
\bigr),
\]
where $M$ is
either:
\begin{longlist}[(ii)]
\item[(i)] equal to $B$; or
\item[(ii)] a bounded martingale orthogonal to $B$.
\end{longlist}
\end{lemma}

\begin{pf}
We prove each claim in turn.
\begin{longlist}[(ii)]
\item[(i)]
For $p,q > 1$, $1/p + 1/q = 1$, we have
\begin{eqnarray*}
&&\bigl\vert\E \bigl[ \bigl(1-K(\tau_n\Delta X_{j,k,n}) \bigr)
(B_{t_{j,k}}-B_{t_j})\mid \mathcal F_{t_j} \bigr]\bigr\vert
\\
&&\qquad\le\E \bigl[\bigl\vert1-K(\tau_n\Delta X_{j,k,n})
\bigr\vert^p\mid\mathcal F_{t_j} \bigr]^{1/p} \E \bigl[
\vert B_{t_{j,k}}-B_{t_j}\vert^q\mid\mathcal
F_{t_j} \bigr]^{1/q},
\end{eqnarray*}
using H\"older's inequality,
\[
=O \bigl(n^{-1/2}
\bigr)\E \bigl[1-K(\tau_n\Delta X_{j,k,n})\mid\mathcal
F_{t_j} \bigr]^{1/p},
\]
using Lemma~\ref{lem:bounds}(ii), and since $K$ takes values in
$[0,1]$,
\[
=O \bigl(n^{-1/2-(1-\alpha\beta)/p} \bigr),
\]
using Lemma~\ref{lem:jmp},
\[
=o \bigl(n^{-(1-\alpha\beta/2)} \bigr),
\]
for small enough $p$. Summing this result, we conclude that
\[
\sum_{j=0}^{\lfloor nt
\rfloor-2m}\E \bigl[ \bigl(1-K(
\tau_n\Delta X_{j,k,n}) \bigr) (B_{t_{j,k}}-B_{t_j})
\mid\mathcal F_{t_j} \bigr] = o \bigl(n^{\alpha
\beta/2} \bigr).
\]

\item[(ii)] Using Lemma~\ref{lem:lvy}, for fixed $k$ and $n$, we can write
\[
\Delta X_{j,k,n} = \overline\xi_{t_{j,k}}^{(j)} +
Y_{1}^{(j)} + Y_{2}^{(j)},
\]
for a L\'evy approximation $\overline
\xi_{t}^{(j)}$, and error terms $Y_{1}^{(j)}$, $Y_{2}^{(j)}$.
We can then write
\begin{eqnarray*}
&&\sum_{j=0}^{\lfloor nt \rfloor- 2m}\E \bigl[ \bigl(1-K(
\tau_n\Delta X_{j,k,n}) \bigr) (M_{t_{j,k}}-M_{t_j})
\mid\mathcal F_{t_j} \bigr]
\\
&&\qquad= %
\sum
_{j=0}^{\lfloor nt \rfloor-
2m}\E \bigl[ \bigl(K \bigl(\tau_n
\overline\xi_{t_{j,k}}^{(j)} \bigr)-K(\tau_n\Delta
X_{j,k,n}) \bigr) (M_{t_{j,k}}-M_{t_j})\mid\mathcal
F_{t_j} \bigr]
\\
&&\qquad\quad{}+ \sum_{j=0}^{\lfloor nt \rfloor- 2m}\E \bigl[ \bigl(1-K
\bigl( \tau_n\overline\xi_{t_{j,k}}^{(j)} \bigr) \bigr)
(M_{t_{j,k}}-M_{t_j})\mid\mathcal F_{t_j} \bigr],
\end{eqnarray*}
where we will bound separately the two sums on the right-hand
side.

For the first sum, we have
\begin{eqnarray*}
&&\sum_{j=0}^{\lfloor nt \rfloor- 2m}\E \bigl[ %
\bigl(K \bigl(\tau_n\overline
\xi_{t_{j,k}}^{(j)} \bigr)-K(\tau_n\Delta
X_{j,k,n}) \bigr) (M_{t_{j,k}}-M_{t_j})\mid\mathcal
F_{t_j} \bigr]  %
\\
&&\qquad= %
O(1)\sum
_{j=0}^{\lfloor nt \rfloor-
2m} \E \bigl[ \bigl(1 \wedge\bigl\vert\tau
_n \bigl(Y_{1}^{(j)} + Y_{2}^{(j)}
\bigr)\bigr\vert \bigr)\vert M_{t_{j,k}}-M_{t_j}\vert\mid\mathcal
F_{t_j} \bigr]
\end{eqnarray*}
since $K(x+y) = K(x) + O(1\wedge\vert y\vert)$,
\begin{eqnarray*}
&&= O(1)\sum_{j=0}^{\lfloor nt \rfloor- 2m}
\bigl( %
\E \bigl[\bigl\vert
\tau_n Y_{1}^{(j)}\bigr\vert ^2\mid
\mathcal F_{t_j} \bigr]^{1/2}\E \bigl[(M_{t_{j,k}}-M_{t_j})^2
\mid \mathcal F_{t_j} \bigr]^{1/2}
\\
&&\quad{}+ \E \bigl[1 \wedge \bigl\vert\tau_n Y_{2}^{(j)}
\bigr\vert\mid\mathcal F_{t_j} \bigr] \bigr), %
\end{eqnarray*}
by Cauchy--Schwarz, and since $M$ is bounded,
\[
=o_p
\bigl(n^{-1/2} \bigr)\sum_{j=0}^{\lfloor nt \rfloor- 2m}
\E \bigl[(M_{t_{j,k}}-M_{t_j})^2 \bigr]^{1/2}
+ o \bigl(n^{\alpha\beta/2} \bigr),
\]
using Lemma~\ref{lem:lvy},
\[=o_p(1) \Biggl(\sum
_{j=0}^{\lfloor nt \rfloor- 2m}\E \bigl[(M_{t_{j,k}}-M_{t_j})^2
\bigr] \Biggr)^{1/2} + o \bigl(n^{\alpha\beta/2} \bigr),
\]
using Cauchy--Schwarz,
\[
=o_p(1)\E \bigl[(M_1-M_0)^2
\bigr]^{1/2} + o \bigl(n^{\alpha\beta/2} \bigr),
\]
as $M$ is a martingale,
\[
=o_p \bigl(n^{\alpha\beta/2} \bigr),
\]
as $M$ is bounded.

It remains to bound the second sum. Given $\mathcal
F_{t_j}$, $\overline\xi_t^{(j)}$ is a function of the
Brownian motion $B$ and Poisson random measure $\mu$, so we
may apply Theorem III.4.34 of \citet{jacodlimit2003}. We deduce
that
\begin{eqnarray*}
&&K \bigl(\tau_n\overline\xi_{t_{j,k}}^{(j)} \bigr) -
\E \bigl[K \bigl(\tau_n\overline \xi _{t_{j,k}}^{(j)}
\bigr) \mid\mathcal F_{t_j} \bigr]
\\
&&\qquad =\int_{t_j}^{t_{j,k}} G_{s}^{(j)}
\,dB_s + \int_{t_j}^{t_{j,k}}\int
_\R {G_{s}'}^{(j)}(x)
\bigl(\mu(dx,ds)-dx\, ds \bigr),
\end{eqnarray*}
for a predictable process $G_{s}^{(j)}$, and predictable
function ${G_{s}'}^{(j)}(x)$. Likewise, by their Lemma
III.4.24, we have
\[
M_{t}-M_{0} = \int_{0}^{t}
G_s''(x) \bigl(\mu(dx,ds)-dx \,ds \bigr)+
\overline M_t,
\]
for a predictable function $G_s''(x)$, and a martingale $\overline
M_t$
orthogonal to $B$ and $\mu$.

Now, as $K$ is bounded, so is ${G'_s}^{(j)}$; furthermore, by
considering the quadratic variation, we have
%
\begin{eqnarray}\label{eq:qva}
&&
\E \biggl[\int_{t_j}^{t_{j,k}}\int
_\R{G_{s}'}^{(j)}(x)^2
\,dx \,ds \Big\mid\mathcal F_{t_j} \biggr]
\nonumber\\
&&\qquad\le\Var \bigl[K \bigl(\tau_n\overline
\xi_{t_{j,k}}^{(j)} \bigr)\mid \mathcal F_{t_j} \bigr]
\nonumber
\\[-8pt]
\\[-8pt]
\nonumber
&&\qquad= O(1) \bigl(\Var \bigl[K(\tau_n \Delta_{j,k,n})
\mid\mathcal F_{t_j} \bigr] + \E \bigl[1 \wedge\bigl\vert\tau_n
\bigl(Y_1^{(j)} + Y_2^{(j)} \bigr)
\bigr\vert^2 \mid \mathcal F_{t_j} \bigr] \bigr)
\\
& &\qquad = O \bigl(n^{-(1-\alpha\beta)} \bigr),\nonumber
\end{eqnarray}
using Lemmas~\ref{lem:lvy} and \ref{lem:ev}.
We likewise have
%
\begin{eqnarray}\label{eq:qvb}
\E \biggl[\int_0^{1}\int
_\R G_s''(x)^2
\,dx \,ds \biggr] &\le& \E \bigl[(M_t-M_0)^2
\bigr]
\nonumber
\\[-8pt]
\\[-8pt]
\nonumber
 &=&O(1),
\end{eqnarray}
as $M$ is bounded.

Setting $\varepsilon_n = n^{-\alpha\beta/4}$, we thus obtain
\begin{eqnarray*}
&&\E \bigl[ \bigl(1-K \bigl(\tau_n\overline\xi_{t_{j,k}}^{(j)}
\bigr) \bigr) (M_{t_{j,k}}-M_{t_j}) \mid\mathcal F_{t_j}
\bigr]
\\
&&\qquad=\E \biggl[\int_{t_j}^{t_{j,k}}\int
_\R {G_{s}'}^{(j)}(x)G_s''(x)
\,dx \,ds \Big\mid\mathcal F_{t_j} \biggr],
\end{eqnarray*}
applying It\^o's lemma,
\begin{eqnarray*}
&\le &
\E \biggl[\int_{t_j}^{t_{j,k}}\int
_\R {G_{s}'}^{(j)}(x)^2
\,dx \,ds \Big\mid\mathcal F_{t_j} \biggr]^{1/2}
\\
&&{}\times\E \biggl[\int_{t_j}^{t_{j,k}}\int
_{\vert G_s''(x)\vert \le
\varepsilon
_n} G_s''(x)^2
\,dx \,ds \Big\mid\mathcal F_{t_j} \biggr]^{1/2}
\\
&&{}+O(1)\E \biggl[\int_{t_j}^{t_{j,k}}\int
_{\vert G_s''(x)\vert >
\varepsilon_n} \bigl\vert G_s''(x)
\bigr\vert \,dx \,ds \Big\mid\mathcal F_{t_j} \biggr],
 \end{eqnarray*}
using Cauchy--Schwarz, and since ${G_{s}'}^{(j)}$
is bounded,
\begin{eqnarray*}
&=& %
 O
\bigl(n^{-(1-\alpha\beta)/2} \bigr)\E \biggl[\int_{t_j}^{t_{j,k}}
\int_{\vert
G_s''(x)\vert
\le\varepsilon_n} G_s''(x)^2
\,dx \,ds \Big\mid\mathcal F_{t_j} \biggr]^{1/2}
\\
&&{}+ O \bigl(n^{\alpha\beta/4} \bigr)\E \biggl[\int_{t_j}^{t_{j,k}}
\int_{\vert
G_s''(x)\vert
> \varepsilon_n} G_s''(x)^2
\,dx \,ds \Big\mid\mathcal F_{t_j} \biggr],
 \end{eqnarray*}
using \eqref{eq:qva}. We thus have
\begin{eqnarray*}
&&\sum_{j=0}^{\lfloor nt \rfloor- 2m}\E \bigl[ \bigl(1-K \bigl(
\tau_n\overline \xi _{t_{j,k}}^{(j)} \bigr) \bigr)
(M_{t_{j,k}}-M_{t_j}) \mid\mathcal F_{t_j} \bigr]
\\
&&\qquad= %
O_p
\bigl(n^{\alpha\beta/2} \bigr)\E \biggl[\int_{0}^{1}
\int_{\vert G_s''(x)\vert
\le\varepsilon_n} G_s''(x)^2
\,dx \,ds \biggr]^{1/2}
\\
&&\qquad\quad{}+ O_p \bigl(n^{\alpha\beta/4} \bigr)\E \biggl[\int
_{0}^{1}\int_{\vert G_s''(x)\vert
> \varepsilon_n}
G_s''(x)^2 \,dx \,ds \biggr],
\end{eqnarray*}
using Cauchy--Schwarz again,
\[=o_p \bigl(n^{\alpha\beta/2}
\bigr),
\]
using \eqref{eq:qvb}. \quad\qed
\end{longlist}
\noqed\end{pf}

We now prove a limit theorem for our jump counts $\widehat A'_n(\tau
)$.

\begin{lemma}
\label{lem:lim}
In the setting of Theorem~\ref{thm:lim}, for $l = 0, 1$, set
\[
\eta_{n,l} = \tau_{n,l}^{\beta/2} \bigl(
\tau_{n,l}^{-\beta}\widehat A_n'(
\tau_{n,l}) - K_{\beta}\overline\Gamma_1 \bigr),\qquad
\tau_{n,l} = \rho^l \tau_n.
\]
Then the random vector
\[
\eta_n \overset{\mathrm{sd}}\to\overline\Gamma_1^{1/2}
\widetilde Z,
\]
where the random variable $\widetilde Z$ is defined, on a suitable extension
of the probability space $(\Omega, \mathcal F, \mathbb P)$, to
satisfy
\[
\widetilde Z \mid\mathcal F \sim N \biggl(0, %
\pmatrix{ \overline
K_{\beta,1}&\overline K_{\beta,\rho}\vspace *{2pt}
\cr
\overline
K_{\beta,\rho}&\overline K_{\beta,1}} %
 \biggr).
\]
\end{lemma}

\begin{pf}
We first make a localisation argument, allowing us to work under
Assumption~\ref{ass:bd}. Since we wish to bound both $c_t$ and its
characteristics, we will localise explicitly. For $k=1,2,\ldots,$
let $\phi_k$ be a smooth bounded function with bounded
derivatives, equal to the identity on $[-k,k]$.

Without loss of generality, we may assume that the stopping times
$T_k$ also localise the processes $b_t, b^c_t, H_t, H'_t,
\gamma^\pm_t, \int_\R1 \wedge\vert\delta_t(x)\vert^{\upsilon_1}
\,dx$ and
$\int_\R1 \wedge\delta^c_t(x)^2 \,dx$. We can then write
\begin{eqnarray*}
X_t^{(k)} &=&%
\int_0^t b_{s\wedge T_k} \,ds + \int
_0^t \phi_k \bigl(c_s^{(k)}
\bigr) \,dB_s + \int_0^t
\gamma_{s\wedge T_k}^+ \,dL^+_s + \int_0^t
\gamma_{s\wedge T_k}^- \,dL^-_s
\\
&&{}+ \int_0^t \int_\R
\delta_{s\wedge T_k}(x) \mu(dx, ds),
\end{eqnarray*}
where
\begin{eqnarray*}
c_t^{(k)} &=& c_0 + \int_0^{t}
b^c_{s\wedge T_k} \,ds + \int_0^{t}
H_{s\wedge T_k} \,dB_s + \int_0^{t}
H'_{s\wedge T_k} \,dB_s'
\\
&&{}+ \int_0^{t} \int_\R
\delta^c_{s\wedge T_k}(x) \bigl(\mu(dx,ds) - 1_{\vert\delta^c_{s\wedge T_k}(x)\vert < 1} \,dx \,ds
\bigr).
\end{eqnarray*}
We
note that $X^{(k)} = X$ eventually almost-surely, so it suffices
to prove our result instead for the processes $X^{(k)}$; an
application of It\^o's lemma shows that these processes satisfy
Assumption~\ref{ass:bd}.

We next define random variables
\begin{eqnarray*}
\zeta_{j,n,l}& =& %
\tau_{n,l}^{-\beta/2}\E \Biggl[\sum_{j'=(j-2m+1)\vee
0}^{j \wedge(n-2m-1)}
\widehat a_{j'}(\tau_{n,l})\Big\mid\mathcal F_{t_{j+1}}
\Biggr]
\\
&&{}-\tau_{n,l}^{-\beta/2}\E \Biggl[\sum
_{j'=(j-2m+1)\vee0}^{j \wedge
(n-2m)-1} \widehat a_{j'}(
\tau_{n,l})\Big\mid\mathcal F_{t_j} \Biggr]
\\
&&{}-\tau_{n,l}^{\beta/2}K_{\beta} n^{-1}\overline
\gamma_{t_j}1_{j <
n-2m},
\end{eqnarray*}
so we may write
\[
\eta_{n,l} = \sum_{j=0}^{n-1}
\zeta_{j,n,l} + \psi_{n,l},
\]
where the term
\[
\psi_{n,l} = - \tau_{n,l}^{\beta/2}K_\beta\int
_0^1 (\overline\gamma_s - \overline
\gamma_{\lfloor ns \rfloor/n}1_{s < 1-2m/n}) \,ds.
\]
Since
\begin{eqnarray*}
\E[\psi_{n,l}] &=& O \bigl(n^{\alpha\beta/2} \bigr) \biggl(\int
_0^{1-2m/n} \E\bigl[\vert\overline\gamma_s -
\overline\gamma _{\lfloor ns \rfloor/n}\vert\bigr] \,ds + \int_{1-2m/n}^1
\E\bigl[\vert \overline \gamma_s\vert\bigr] \,ds \biggr)
\\
&= &O \bigl(n^{\alpha\beta/2} \bigr) \biggl(\sum_{* \in\{+,-\}}
\int_0^1 \E \bigl[\bigl\vert\gamma_s^*
- \gamma _{\lfloor ns \rfloor/n}^*\bigr\vert^2 \bigr]^{(1 \wedge\beta)/2} \,ds +
n^{-1} \biggr),
\end{eqnarray*}
as the function $x \mapsto\vert x\vert^\beta$ is $(1 \wedge
\beta)$-Lipschitz, and $ \overline\gamma_s$ is bounded,
\[ = O
\bigl(n^{-(1\wedge\beta-\alpha\beta)/2} \bigr),
\]
as the $\gamma_s^\pm$ are smooth in quadratic
mean,
\[ = o(1),
\]
we deduce that
\[
\eta_{n,l} = \sum_{j=0}^{n-1}
\zeta_{j,n,l} + o_p(1).
\]

The desired result then follows from Theorem~2.2.15 of
\citet{jacoddiscretization2011}, provided that for $t \in[0,1]$,
$l = 0,1$:
\begin{enumerate}[(iii)]
\item[(i)]
$\sum_{j=0}^{n-1} \vert\E[\zeta_{j,n,l}\mid
\mathcal F_{t_j}]\vert \overset{p}\to0$;
\item[(ii)]$\sum_{j=0}^{\lfloor nt \rfloor-1} \Var[\zeta_{j,n,l}\mid
\mathcal
F_{t_j}] \overset{p}\to
\overline K_{\beta,1}\overline\Gamma_t$;
\item[(iii)]$\sum_{j=0}^{\lfloor nt \rfloor-1} \Cov[\zeta_{j,n,0},\zeta
_{j,n,1}\mid
\mathcal F_{t_j}] \overset{p}\to
\overline K_{\beta,\rho}\overline\Gamma_t$;
\item[(iv)]$\sum_{j=0}^{n-1} \E[\vert\zeta_{j,n,l}\vert^p\mid
\mathcal F_{t_j}] \overset{p}\to0$, for
some $p > 2$; and
\item[(v)]
$\sum_{j=0}^{\lfloor nt \rfloor-1}
\E[\zeta_{j,n,l}(M_{t_{j+1}}-M_{t_j})\mid\mathcal F_{t_j}] \overset
{p}\to
0$, where $M$ is either:
\begin{enumerate}[(a)]
\item[(a)] equal to $B$; or
\item[(b)] a bounded martingale orthogonal to $B$.
\end{enumerate}
\end{enumerate}

We now prove each claim in turn.
\begin{longlist}
\item[(i)] From Lemma~\ref{lem:em}, we have that for
$j= 0, \ldots, n-2m-1$,
\begin{eqnarray*}
\E[\zeta_{j,n,l}\mid\mathcal F_{t_j}] &= &\tau_{n,l}^{-\beta/2}
\E \bigl[\widehat a_{j,n}(\tau_{n,l})\mid \mathcal
F_{t_j} \bigr] - \tau_{n,l}^{\beta/2}K_{\beta}
n^{-1}\overline\gamma_{t_j}
\\
&= &o \bigl(n^{-1} \bigr).
\end{eqnarray*}
From the definitions, we also have that for
$j= n-2m, \ldots, n-1$,
\[
\E[\zeta_{j,n,l}\mid\mathcal F_{t_j}] = 0.
\]
We conclude that
\[
\sum_{j=0}^{n-1}\bigl\vert\E[
\zeta_{j,n,l}\mid\mathcal F_{t_j}]\bigr\vert = o(1).
\]
\item[(ii)] From Lemma~\ref{lem:ev}, we have that for
$j= 2m-1, \ldots, n-2m-1$, and terms $Y_{j,n,l}$ satisfying
$\E[\vert Y_{j,n,l}\vert] = o(n^{-1})$,
\begin{eqnarray*}
\Var[\zeta_{j,n,l} \mid\mathcal F_{t_j}] &=& %
\tau_{n,l}^{-\beta}\sum
_{k',k''=1}^m w_{k'} w_{k''}
\\
&&{}\times \sum
_{j'=j-2k'+1}^j \sum
_{j''=j-2k''+1}^j\Cov \bigl[ %
\E \bigl[K(
\tau_{n,l} \Delta X_{j',k',n}) \mid\mathcal F_{t_{j+1}}
\bigr]
\\
&&{}\times \E \bigl[K(\tau_{n,l} \Delta X_{j'',k'',n}) \mid\mathcal
F_{t_{j+1}} \bigr] \mid\mathcal F_{t_j} \bigr]
\\
&=&n^{-1} \Biggl(2\sum_{k=1}^mkw_k
\Biggr)^2 \overline K_{\beta,1}\overline\gamma_{t_j}
+ Y_{j,n,l}
\\
&= &n^{-1} \overline K_{\beta,1}\overline\gamma_{t_j} +
Y_{j,n,l},
\end{eqnarray*}
since
%
\begin{eqnarray}
\label{eq:wks}
2\sum_{k=1}^m k
w_k &=& 1-\sum_{k=0}^m
(-1)^k \pmatrix{m \cr k}\nonumber
\\
&=& 1-(1-1)^m
\\
&=&1.\nonumber
\end{eqnarray}
For $j= 0, \ldots, 2m-2$ or $j= n-2m, \ldots, n-1$, by a similar
argument, we have the same result for terms $Y_{j,n,l}$
satisfying $\E[\vert Y_{j,n,l}\vert] = O(n^{-1})$. We deduce that
\begin{eqnarray*}
\sum_{j=0}^{\lfloor nt \rfloor-1}\Var[\zeta_{j,n,l}
\mid\mathcal F_{t_j}] &=& n^{-1} \overline K_{\beta,1}
\sum_{j=0}^{\lfloor nt \rfloor-1}\overline
\gamma_{t_j} + \sum_{j=0}^{\lfloor nt \rfloor-1}Y_{j,n,l}
\\
&=&\overline K_{\beta,1}\overline\Gamma_t +
o_p(1) + O_p \Biggl( \sum_{j=0}^{\lfloor nt \rfloor-1}
\E\bigl[\vert Y_{j,n,l}\vert\bigr] \Biggr)
\\
&=& \overline K_{\beta,1}\overline\Gamma_t +
o_p(1).
\end{eqnarray*}
\item[(iii)] The result follows similar to part (ii).
\item[(iv)] Since $\zeta_{j,n,l} = O(n^{-\alpha\beta/2})$, the result is trivial
for large enough $p$.
\item[(v)] In either case (a) or (b), we have
\begin{eqnarray*}
&&\sum_{j=0}^{\lfloor nt \rfloor
-1}\E \bigl[
\zeta_{j,n,l}(M_{t_{j+1}}-M_{t_j})\mid\mathcal
F_{t_j} \bigr]
\\
&&\qquad=\tau_{n,l}^{-\beta/2}\sum_{j=0}^{\lfloor nt \rfloor
-1}
\sum_{j'=(j-2m+1) \vee0}^{j \wedge(n-2m-1)}\E \bigl[\widehat
a_{j'}(\tau_{n,l}) (M_{t_{j+1}}-M_{t_j})\mid
\mathcal F_{t_j} \bigr],
\end{eqnarray*}
since $M$ is a martingale,
\[
=\tau_{n,l}^{-\beta/2}\sum
_{j=0}^{\lfloor nt
\rfloor-2m}\E \bigl[\widehat
a_j( \tau_{n,l}) (M_{t_{j,m}}-M_{t_j})\mid
\mathcal F_{t_j} \bigr] + o_p(1),
\]
since $\widehat a_{j,n}(\tau)$ is bounded,
\begin{eqnarray*}
&=&\tau_{n,l}^{-\beta/2}\sum
_{j=0}^{\lfloor nt
\rfloor-2m}\sum_{k=1}^mw_k
\\
&&{}\times\E \bigl[ \bigl(1-K(\tau_{n,l}\Delta X_{j,k,n}) \bigr)
(M_{t_{j,k}}-M_{t_j})\mid \mathcal F_{t_j} \bigr] +
o_p(1),
\end{eqnarray*}
from the definition of $\widehat a_{j,n}(\tau)$,
\[
=o_p(1),
\]
using Lemma~\ref{lem:cor}.\quad\qed
\end{longlist}
\noqed\end{pf}

Finally, we can prove a limit theorem for $\widehat\beta_n$.

\begin{pf*}{Proof of Theorem~\ref{thm:lim}}
We begin by defining the variables
\[
\widehat\beta_n' = \log_\rho \biggl(
\frac{A'_n(\rho\tau_n)}{A'_n(\tau_n)} \biggr),\qquad \widehat U_n'(\beta) =
\frac{\tau_n^{\widehat\beta_n'/2}(\widehat
\beta_n' -
\beta)}{\widehat\sigma_{\rho,n}}.
\]
From Lemma~\ref{lem:lim}, on the event $\overline\Gamma_1 > 0$, we
have that
\[
A_n'(\tau_{n,l}) \overset{p}\to
\tau_{n,l}K_\beta\overline \Gamma_1, \qquad \widehat
\beta_n' \overset{p}\to\beta.
\]
Hence,
with probability tending to one,
\[
\widehat\beta_n' = \widehat\beta_n,\qquad
\widehat U_n'(\beta) = \widehat U_n(
\beta).
\]
It thus suffices to prove limit theorems
for the quantities $\widehat\beta_n'$ and $\widehat U_n'(\beta)$.

Next, we note we may equivalently define
$\sigma^2_{\beta,\rho}$ by
\[
\sigma_{\beta,\rho}^2 = \bigl( \bigl(1+\rho^{-\beta} \bigr)
\overline K_{\beta,1}-2\rho^{-\beta/2}\overline K_{\beta,\rho}
\bigr)/ \log(\rho)^2K_\beta^2\overline
\Gamma_1.
\]
Again using Lemma~\ref{lem:lim}, on the event $\overline\Gamma_1 >
0$, we also have\vspace*{-1pt}
\begin{eqnarray*}
\widehat\beta_n' &=& \log_\rho \bigl(
\widehat A_n'(\rho\tau_n)/\widehat
A_n'(\tau_n) \bigr)
\\[-1pt]
&=& \beta+ \log_\rho \bigl(1+\eta_{n,1}/(\rho
\tau_{n})^{\beta/2}K_\beta \overline\Gamma_1
\bigr) -\log_\rho \bigl(1+\eta_{n,0}/\tau_{n}^{\beta/2}K_\beta
\overline\Gamma _1 \bigr)
\\[-1pt]
&=& \beta+ \bigl(\rho^{-\beta/2}\eta_{n,1} - \eta_{n,0}
\bigr)/\log(\rho)\tau_{n}^{\beta/2}K_\beta\overline
\Gamma_1 + o_p \bigl(n^{-\alpha\beta/2} \bigr);\vspace*{-1pt}
\end{eqnarray*}
we deduce that\vspace*{-1pt}
\[
\tau_n^{\beta/2} \bigl(\widehat\beta_n'
- \beta \bigr) \overset{\mathrm{sd}}\to \sigma_{\beta,\rho} Z.
\]
Similarly, we have that $\widehat\sigma_{\rho,n}^2
\overset{p}\to\sigma_{\beta,\rho}^2$. Using equation~(2.2.5) of
\citet{jacoddiscretization2011}, we thus obtain\vspace*{-1pt}
\[
\bigl(\tau_n^{\beta/2} \bigl(\widehat\beta_n'
- \beta \bigr), \widehat \sigma_{\rho,n}^2 \bigr) \overset{\mathrm{sd}}
\to \bigl(\sigma_{\beta,\rho} Z, \sigma_{\beta,\rho}^2 \bigr).\vspace*{-1pt}
\]
By continuous mapping, we deduce that\vspace*{-1pt}
\[
\widehat U_n'(\beta) \overset{\mathrm{sd}}\to Z.\vspace*{-1pt}
\]
\upqed\end{pf*}

We have thus proved Theorem~\ref{thm:lim}; we note that Corollary~\ref{cor:lim}
then follows directly.


\section*{Acknowledgments}
We thank the anonymous referees for
  their valuable suggestions and comments. All research data was randomly generated
  using software given in \citet{bullsoftware2015}.\vspace*{-9pt}

\begin{supplement}[id=suppA]
\stitle{Supplement to ``Near-optimal estimation of jump activity in
semimartingales''}
\slink[doi]{10.1214/15-AOS1349SUPP} 
\sdatatype{.pdf}
\sfilename{aos1349\_supp.pdf}
\sdescription{We provide proofs of our technical results.}
\end{supplement}\vspace*{-9pt}


\begin{thebibliography}{29}

\bibitem[\protect\citeauthoryear{A{\"{\i}}t-Sahalia and
Jacod}{2009}]{ait-sahaliaestimating2009}
%
\begin{barticle}[mr]
\bauthor{\bsnm{A{\"{\i}}t-Sahalia},~\bfnm{Yacine}\binits{Y.}} \AND
\bauthor{\bsnm{Jacod},~\bfnm{Jean}\binits{J.}}
(\byear{2009}).
\btitle{Estimating the degree of activity of jumps in high frequency data}.
\bjournal{Ann. Statist.}
\bvolume{37}
\bpages{2202--2244}.
\bid{doi={10.1214/08-AOS640}, issn={0090-5364}, mr={2543690}}
\end{barticle}
%

\bptok{imsref}%
\endbibitem

\bibitem[\protect\citeauthoryear{A{\"{\i}}t-Sahalia and
Jacod}{2011}]{ait-sahaliatesting2012}
%
\begin{barticle}[mr]
\bauthor{\bsnm{A{\"{\i}}t-Sahalia},~\bfnm{Yacine}\binits{Y.}} \AND
\bauthor{\bsnm{Jacod},~\bfnm{Jean}\binits{J.}}
(\byear{2011}).
\btitle{Testing whether jumps have finite or infinite activity}.
\bjournal{Ann. Statist.}
\bvolume{39}
\bpages{1689--1719}.
\bid{doi={10.1214/11-AOS873}, issn={0090-5364}, mr={2850217}}
\end{barticle}
%

\bptok{imsref}%
\endbibitem

\bibitem[\protect\citeauthoryear{A{\"{\i}}t-Sahalia and
Jacod}{2012}]{ait-sahaliaidentifying2012}
%
\begin{barticle}[mr]
\bauthor{\bsnm{A{\"{\i}}t-Sahalia},~\bfnm{Yacine}\binits{Y.}} \AND
\bauthor{\bsnm{Jacod},~\bfnm{Jean}\binits{J.}}
(\byear{2012}).
\btitle{Identifying the successive {B}lumenthal--{G}etoor indices of a
discretely observed process}.
\bjournal{Ann. Statist.}
\bvolume{40}
\bpages{1430--1464}.
\bid{doi={10.1214/12-AOS976}, issn={0090-5364}, mr={3015031}}
\end{barticle}
%

\bptok{imsref}%
\endbibitem

\bibitem[\protect\citeauthoryear{Barndorff-Nielsen and
Shephard}{2004}]{barndorff-nielsenpower2004}
%
\begin{barticle}[author]
\bauthor{\bsnm{Barndorff-Nielsen},~\bfnm{Ole~E.}\binits{O.~E.}}
\AND
\bauthor{\bsnm{Shephard},~\bfnm{Neil}\binits{N.}}
(\byear{2004}).
\btitle{Power and bipower variation with stochastic volatility and jumps}.
\bjournal{Journal of Financial Econometrics}
\bvolume{2}
\bpages{1--37}.
\end{barticle}
%

\bptok{imsref}%
\endbibitem

\bibitem[\protect\citeauthoryear{Barndorff-Nielsen
et~al.}{2006}]{barndorff-nielsenlimit2006}
%
\begin{barticle}[mr]
\bauthor{\bsnm{Barndorff-Nielsen},~\bfnm{Ole~E.}\binits{O.~E.}},
\bauthor{\bsnm{Graversen},~\bfnm{Svend~Erik}\binits{S.~E.}},
\bauthor{\bsnm{Jacod},~\bfnm{Jean}\binits{J.}} \AND
\bauthor{\bsnm{Shephard},~\bfnm{Neil}\binits{N.}}
(\byear{2006}).
\btitle{Limit theorems for bipower variation in financial econometrics}.
\bjournal{Econometric Theory}
\bvolume{22}
\bpages{677--719}.
\bid{doi={10.1017/S0266466606060324}, issn={0266-4666}, mr={2283032}}
\end{barticle}
%

\bptok{imsref}%
\endbibitem

\bibitem[\protect\citeauthoryear{Blumenthal and
Getoor}{1961}]{blumenthalsample1961}
%
\begin{barticle}[mr]
\bauthor{\bsnm{Blumenthal},~\bfnm{R.~M.}\binits{R.~M.}} \AND
\bauthor{\bsnm{Getoor},~\bfnm{R.~K.}\binits{R.~K.}}
(\byear{1961}).
\btitle{Sample functions of stochastic processes with stationary
independent increments}.
\bjournal{J. Math. Mech.}
\bvolume{10}
\bpages{493--516}.
\bid{mr={0123362}}
\end{barticle}\vadjust{\eject}
%

\bptok{imsref}%
\endbibitem

\bibitem[\protect\citeauthoryear{Bull}{2014}]{bullestimating2013}
%
\begin{barticle}[mr]
\bauthor{\bsnm{Bull},~\bfnm{Adam~D.}\binits{A.~D.}}
(\byear{2014}).
\btitle{Estimating time-changes in noisy L\'evy models}.
\bjournal{Ann. Statist.}
\bvolume{42}
\bpages{2026--2057}.
\bid{doi={10.1214/14-AOS1250}, issn={0090-5364}, mr={3262476}}\vadjust{\goodbreak}
\end{barticle}
%

\bptok{imsref}%
\endbibitem

\bibitem[\protect\citeauthoryear{Bull}{2015a}]{bullsoftware2015}
%
\begin{bmisc}[author]
\bauthor{\bsnm{Bull},~\bfnm{Adam~D.}\binits{A.~D.}}
(\byear{2015}a).
\bhowpublished{Software for ``{N}ear-optimal estimation of jump
activity in semimartingales.''
Available at \url{https://www.repository.cam.ac.uk/handle/1810/248959}}.
\end{bmisc}
%
\bptok{imsref}%
\endbibitem

\bibitem[\protect\citeauthoryear{Bull}{2015b}]{bullsupplement2015}
%
\begin{bmisc}[author]
{\bauthor{\bsnm{Bull},~\binits{A. D.}}}
(\byear{2015b}).
\bhowpublished{Supplement to ``Near-optimal estimation of jump
activity in semimartingales.''
DOI:\doiurl{10.1214/15-AOS1349SUPP}}.
\bptok{imsref}%
\end{bmisc}
%
\endbibitem

\bibitem[\protect\citeauthoryear{Cont and
Mancini}{2011}]{contnonparametric2011}
%
\begin{barticle}[mr]
\bauthor{\bsnm{Cont},~\bfnm{Rama}\binits{R.}} \AND
\bauthor{\bsnm{Mancini},~\bfnm{Cecilia}\binits{C.}}
(\byear{2011}).
\btitle{Nonparametric tests for pathwise properties of semimartingales}.
\bjournal{Bernoulli}
\bvolume{17}
\bpages{781--813}.
\bid{doi={10.3150/10-BEJ293}, issn={1350-7265}, mr={2787615}}
\end{barticle}
%

\bptok{imsref}%
\endbibitem

\bibitem[\protect\citeauthoryear{Cont and Tankov}{2004}]{contfinancial2004}
%
\begin{bbook}[mr]
\bauthor{\bsnm{Cont},~\bfnm{Rama}\binits{R.}} \AND
\bauthor{\bsnm{Tankov},~\bfnm{Peter}\binits{P.}}
(\byear{2004}).
\btitle{Financial Modelling with Jump Processes}.
\bpublisher{Chapman \& Hall/CRC},
\blocation{Boca Raton, FL}.
\bid{mr={2042661}}
\end{bbook}
%

\bptok{imsref}%
\endbibitem

\bibitem[\protect\citeauthoryear{Jacod and
Protter}{2012}]{jacoddiscretization2011}
%
\begin{bbook}[mr]
\bauthor{\bsnm{Jacod},~\bfnm{Jean}\binits{J.}} \AND
\bauthor{\bsnm{Protter},~\bfnm{Philip}\binits{P.}}
(\byear{2012}).
\btitle{Discretization of Processes}.
\bseries{Stochastic Modelling and Applied Probability}
\bvolume{67}.
\bpublisher{Springer},
\blocation{Heidelberg}.
\bid{doi={10.1007/978-3-642-24127-7}, mr={2859096}}
\end{bbook}
%

\bptok{imsref}%
\endbibitem

\bibitem[\protect\citeauthoryear{Jacod and Reiss}{2014}]{jacodremark2012}
%
\begin{barticle}[mr]
\bauthor{\bsnm{Jacod},~\bfnm{Jean}\binits{J.}} \AND
\bauthor{\bsnm{Reiss},~\bfnm{Markus}\binits{M.}}
(\byear{2014}).
\btitle{A remark on the rates of convergence for integrated volatility
estimation in the presence of jumps}.
\bjournal{Ann. Statist.}
\bvolume{42}
\bpages{1131--1144}.
\bid{doi={10.1214/13-AOS1179}, issn={0090-5364}, mr={3224283}}
\end{barticle}
%

\bptok{imsref}%
\endbibitem

\bibitem[\protect\citeauthoryear{Jacod and Shiryaev}{2003}]{jacodlimit2003}
%
\begin{bbook}[mr]
\bauthor{\bsnm{Jacod},~\bfnm{Jean}\binits{J.}} \AND
\bauthor{\bsnm{Shiryaev},~\bfnm{Albert~N.}\binits{A.~N.}}
(\byear{2003}).
\btitle{Limit Theorems for Stochastic Processes},
\bedition{2nd} ed.
\bseries{Grundlehren der Mathematischen Wissenschaften}
\bvolume{288}.
\bpublisher{Springer},
\blocation{Berlin}.
\bid{doi={10.1007/978-3-662-05265-5}, mr={1943877}}
\end{bbook}
%

\bptok{imsref}%
\endbibitem

\bibitem[\protect\citeauthoryear{Jacod and
Todorov}{2014}]{jacodefficient2014}
%
\begin{barticle}[mr]
\bauthor{\bsnm{Jacod},~\bfnm{Jean}\binits{J.}} \AND
\bauthor{\bsnm{Todorov},~\bfnm{Viktor}\binits{V.}}
(\byear{2014}).
\btitle{Efficient estimation of integrated volatility in presence of
infinite variation jumps}.
\bjournal{Ann. Statist.}
\bvolume{42}
\bpages{1029--1069}.
\bid{doi={10.1214/14-AOS1213}, issn={0090-5364}, mr={3210995}}
\end{barticle}
%

\bptok{imsref}%
\endbibitem

\bibitem[\protect\citeauthoryear{Jing, Kong and
Liu}{2011}]{jingestimating2011}
%
\begin{barticle}[mr]
\bauthor{\bsnm{Jing},~\bfnm{Bing-Yi}\binits{B.-Y.}},
\bauthor{\bsnm{Kong},~\bfnm{Xin-Bing}\binits{X.-B.}} \AND
\bauthor{\bsnm{Liu},~\bfnm{Zhi}\binits{Z.}}
(\byear{2011}).
\btitle{Estimating the jump activity index under noisy observations
using high-frequency data}.
\bjournal{J. Amer. Statist. Assoc.}
\bvolume{106}
\bpages{558--568}.
\bid{doi={10.1198/jasa.2011.tm10021}, issn={0162-1459}, mr={2847970}}
\end{barticle}
%

\bptok{imsref}%
\endbibitem

\bibitem[\protect\citeauthoryear{Jing et~al.}{2012}]{jingjump2012}
%
\begin{barticle}[mr]
\bauthor{\bsnm{Jing},~\bfnm{Bing-Yi}\binits{B.-Y.}},
\bauthor{\bsnm{Kong},~\bfnm{Xin-Bing}\binits{X.-B.}},
\bauthor{\bsnm{Liu},~\bfnm{Zhi}\binits{Z.}} \AND
\bauthor{\bsnm{Mykland},~\bfnm{Per}\binits{P.}}
(\byear{2012}).
\btitle{On the jump activity index for semimartingales}.
\bjournal{J. Econometrics}
\bvolume{166}
\bpages{213--223}.
\bid{doi={10.1016/j.jeconom.2011.09.036}, issn={0304-4076}, mr={2862961}}
\end{barticle}
%

\bptok{imsref}%
\endbibitem

\bibitem[\protect\citeauthoryear{Lee and Hannig}{2010}]{leedetecting2010}
%
\begin{barticle}[author]
\bauthor{\bsnm{Lee},~\bfnm{Suzanne~S.}\binits{S.~S.}} \AND
\bauthor{\bsnm{Hannig},~\bfnm{Jan}\binits{J.}}
(\byear{2010}).
\btitle{Detecting jumps from {L\'e}vy jump diffusion processes}.
\bjournal{Journal of Financial Economics}
\bvolume{96}
\bpages{271--290}.
\end{barticle}
%

\bptok{imsref}%
\endbibitem

\bibitem[\protect\citeauthoryear{Mancini}{2001}]{mancinidisentangling2001}
%
\begin{barticle}[author]
\bauthor{\bsnm{Mancini},~\bfnm{Cecilia}\binits{C.}}
(\byear{2001}).
\btitle{Disentangling the jumps of the diffusion in a geometric
jumping {B}rownian motion}.
\bjournal{Giornale dell'Istituto Italiano degli Attuari}
\bvolume{64}
\bpages{19--47}.
\end{barticle}
%

\bptok{imsref}%
\endbibitem

\bibitem[\protect\citeauthoryear{Mancini}{2009}]{mancininon2006}
%
\begin{barticle}[mr]
\bauthor{\bsnm{Mancini},~\bfnm{Cecilia}\binits{C.}}
(\byear{2009}).
\btitle{Nonparametric threshold estimation for models with stochastic
diffusion coefficient and jumps}.
\bjournal{Scand. J. Stat.}
\bvolume{36}
\bpages{270--296}.
\bid{doi={10.1111/j.1467-9469.2008.00622.x}, issn={0303-6898}, mr={2528985}}
\end{barticle}
%

\bptok{imsref}%
\endbibitem

\bibitem[\protect\citeauthoryear
{Papapantoleon}{2008}]{papapantoleonintroduction2008}
%
\begin{bmisc}[author]
\bauthor{\bsnm{Papapantoleon},~\bfnm{Antonis}\binits{A.}}
(\byear{2008}).
\bhowpublished{An introduction to {L\'e}vy processes with applications
in finance.
Preprint. Available at \arxivurl{arXiv:0804.0482}}.
\end{bmisc}
%

\bptok{imsref}%
\endbibitem

\bibitem[\protect\citeauthoryear{Podolskij and
Vetter}{2009}]{podolskijestimation2009}
%
\begin{barticle}[mr]
\bauthor{\bsnm{Podolskij},~\bfnm{Mark}\binits{M.}} \AND
\bauthor{\bsnm{Vetter},~\bfnm{Mathias}\binits{M.}}
(\byear{2009}).
\btitle{Estimation of volatility functionals in the simultaneous
presence of microstructure noise and jumps}.
\bjournal{Bernoulli}
\bvolume{15}
\bpages{634--658}.
\bid{doi={10.3150/08-BEJ167}, issn={1350-7265}, mr={2555193}}
\end{barticle}
%

\bptok{imsref}%
\endbibitem

\bibitem[\protect\citeauthoryear{Rei{\ss}}{2013}]{reistesting2013}
%
\begin{barticle}[mr]
\bauthor{\bsnm{Rei{\ss}},~\bfnm{Markus}\binits{M.}}
(\byear{2013}).
\btitle{Testing the characteristics of a L\'evy process}.
\bjournal{Stochastic Process. Appl.}
\bvolume{123}
\bpages{2808--2828}.
\bid{doi={10.1016/j.spa.2013.03.016}, issn={0304-4149}, mr={3054546}}
\end{barticle}
%

\bptok{imsref}%
\endbibitem

\bibitem[\protect\citeauthoryear{Sato}{1999}]{satolevy1999}
%
\begin{bbook}[mr]
\bauthor{\bsnm{Sato},~\bfnm{Ken-iti}\binits{K.-i.}}
(\byear{1999}).
\btitle{L\'evy Processes and Infinitely Divisible Distributions}.
\bseries{Cambridge Studies in Advanced Mathematics}
\bvolume{68}.
\bpublisher{Cambridge Univ. Press},
\blocation{Cambridge}.
\bid{mr={1739520}}
\end{bbook}
%

\bptok{imsref}%
\endbibitem

\bibitem[\protect\citeauthoryear{Todorov and
Tauchen}{2010}]{todorovactivity2010}
%
\begin{barticle}[mr]
\bauthor{\bsnm{Todorov},~\bfnm{Viktor}\binits{V.}} \AND
\bauthor{\bsnm{Tauchen},~\bfnm{George}\binits{G.}}
(\byear{2010}).
\btitle{Activity signature functions for high-frequency data analysis}.
\bjournal{J. Econometrics}
\bvolume{154}
\bpages{125--138}.
\bid{doi={10.1016/j.jeconom.2009.06.009}, issn={0304-4076}, mr={2558956}}
\end{barticle}
%

\bptok{imsref}%
\endbibitem

\bibitem[\protect\citeauthoryear{Todorov and
Tauchen}{2012a}]{todorovinverse2012}
%
\begin{barticle}[mr]
\bauthor{\bsnm{Todorov},~\bfnm{Viktor}\binits{V.}} \AND
\bauthor{\bsnm{Tauchen},~\bfnm{George}\binits{G.}}
(\byear{2012}a).
\btitle{Inverse realized {L}aplace transforms for nonparametric
volatility density estimation in jump-diffusions}.
\bjournal{J. Amer. Statist. Assoc.}
\bvolume{107}
\bpages{622--635}.
\bid{doi={10.1080/01621459.2012.682854}, issn={0162-1459}, mr={2980072}}
\end{barticle}
%

\bptok{imsref}%
\endbibitem

\bibitem[\protect\citeauthoryear{Todorov and
Tauchen}{2012b}]{todorovrealized2012}
%
\begin{barticle}[mr]
\bauthor{\bsnm{Todorov},~\bfnm{Viktor}\binits{V.}} \AND
\bauthor{\bsnm{Tauchen},~\bfnm{George}\binits{G.}}
(\byear{2012}b).
\btitle{The realized {L}aplace transform of volatility}.
\bjournal{Econometrica}
\bvolume{80}
\bpages{1105--1127}.
\bid{doi={10.3982/ECTA9133}, issn={0012-9682}, mr={2963883}}
\end{barticle}
%

\bptok{imsref}%
\endbibitem

\bibitem[\protect\citeauthoryear{Woerner}{2011}]{woerneranalyzing2011}
%
\begin{bincollection}[mr]
\bauthor{\bsnm{Woerner},~\bfnm{Jeannette~H.~C.}\binits{J.~H.~C.}}
(\byear{2011}).
\btitle{Analyzing the fine structure of continuous time stochastic processes}.
In \bbooktitle{Seminar on {S}tochastic {A}nalysis, {R}andom {F}ields
and {A}pplications {VI}}.
\bseries{Progress in Probability}
\bvolume{63}
\bpages{473--492}.
\bpublisher{Birkh\"auser},
\blocation{Basel}.
\bid{doi={10.1007/978-3-0348-0021-1_26}, mr={2857040}}
\end{bincollection}\vadjust{\eject}
%

\bptok{imsref}%
\endbibitem

\bibitem[\protect\citeauthoryear{Zhao and Wu}{2009}]{zhaononparametric2009}
%
\begin{barticle}[mr]
\bauthor{\bsnm{Zhao},~\bfnm{Zhibiao}\binits{Z.}} \AND
\bauthor{\bsnm{Wu},~\bfnm{Wei~Biao}\binits{W.~B.}}
(\byear{2009}).
\btitle{Nonparametric inference of discretely sampled stable L\'evy processes}.
\bjournal{J. Econometrics}
\bvolume{153}
\bpages{83--92}.
\bid{doi={10.1016/j.jeconom.2009.04.007}, issn={0304-4076}, mr={2558496}}
\end{barticle}
%

\bptok{imsref}%
\endbibitem
\end{thebibliography}
%
%



\printaddresses
\end{document}